\documentclass[11pt,a4paper]{article}
\usepackage{authblk}
\usepackage{amsfonts,amsmath,amssymb,mathrsfs,url,esint,comment
}

\usepackage[normalem]{ulem}
\usepackage{hyperref}
\usepackage{cleveref}

\usepackage{color}

\usepackage{graphicx}
\usepackage{epsfig}
\usepackage{float}
\usepackage{pdfsync}
\usepackage{wasysym}
\usepackage{empheq}
\usepackage{textgreek}
\usepackage{url}
\usepackage[left=2cm,right=2cm,top=2cm,bottom=2.5cm]{geometry}
\newcommand{\stkout}[1]{\ifmmode\text{\sout{\ensuremath{#1}}}\else\sout{#1}\fi}

\newcommand{\RN}[1]{\textup{\uppercase\expandafter{\romannumeral#1}}}
\newcommand{\dis}{\displaystyle}

\def\e{{\mathop {\rm e}\nolimits}}

\newcommand\bsout{\bgroup\markoverwith{\textcolor{blue}{\rule[0.5ex]{2pt}{0.4pt}}}\ULon}

\newcommand{\ub}{{\bar u}}

\newcommand{\wto}{\rightharpoonup}
\newcommand{\weaks}{\stackrel{*}{\wto}}

\newcommand{\eps}{\varepsilon}

\newcommand{\R}{\mathbb{R}}
\newcommand{\N}{\mathbb{N}}

\providecommand{\dg}[1]{\textcolor{blue}{#1}}

\def\umax{{u_{\max}}}

\newcommand{\tae}{\tau^{T}}
\newcommand{\sie}{\sigma^{T}}

\providecommand{\pr}[1]{\left(#1\right)} 
\providecommand{\pp}[1]{\left[#1\right]} 

\graphicspath{{../vlasov/},{../Innovaction_e_precorso_2008/},{../../../figure/},{../../../../DIDATTICA/Dispense/figure/}}

\newtheorem{theorem}{Theorem}[section]
\newtheorem{lemma}[theorem]{Lemma}
\newtheorem{proposition}[theorem]{Proposition}
\newtheorem{problem}{Problem}[section]
\newtheorem{corollary}[theorem]{Corollary}
\newtheorem{remark}[theorem]{Remark}
\newtheorem{example}[theorem]{Example}
\newtheorem{definition}[theorem]{Definition}

\def\proof{\noindent{\textbf{Proof. }}}
\def\QED{\hfill {$\square$}\goodbreak \medskip}

\renewcommand{\theta}{\vartheta}

\numberwithin{equation}{section}
\def\fe{for example } 
\def\bep{\begin{pmatrix}} \def\eep{\end{pmatrix}}\newcommand{\be}[1]{\begin{equation}\label{#1}}
\newcommand{\ee}{\end{equation}}
\def\bea{\begin{eqnarray*}}
    
\def\eea{\end{eqnarray*}}    \def\fe{for example } 
\date{}

\title 
{Infinite horizon optimal control of a SIR epidemic\\ under an ICU constraint}
\author[2]{Lorenzo Freddi}
\author[3,4]{Dan Goreac}
\affil[2]{Dipartimento di Scienze Matematiche, Informatiche e Fisiche,  Universit\`a di Udine, via delle Scienze 206, 33100 Udine, Italy,  lorenzo.freddi@uniud.it}
\affil[3]{School of Mathematics and Statistics, Shandong University, Weihai, Weihai 264209, PR China}
\affil[4]{LAMA, Univ Gustave Eiffel, UPEM, Univ Paris Est Creteil, CNRS, F-77447 Marne-la-Valle\'ee, France, dan.goreac@univ-eiffel.fr}

\begin{document}

\maketitle
\abstract{
The aim of this paper is to provide a rigorous mathematical analysis of an optimal control problem of a SIR epidemic on an infinite horizon. 
A state constraint related to intensive care units (ICU) capacity is imposed and the objective functional linearly depends on the state and the control. After preliminary asymptotic and viability analyses, a $\Gamma$-convergence argument is developed to reduce the problem to a finite horizon allowing to use a state constrained version of Pontryagin's theorem to characterize the structure of the optimal controls. 
Illustrating examples and numerical simulations are given according to the available data on Covid-19 epidemic in Italy.
\\

}
\section{Introduction}

The optimal control of epidemics   (\cite{anderson1992infectious,Behncke,hansen2011optimal,Mart}) have attracted the interest 
of researchers since the introduction of the first compartmental epidemic model by Kermack and McKendrick \cite{kermack1927contribution}.
Lots of new contributions have appeared in recent years, particularly after the start of the COVID-19 pandemic -- see \fe\ \cite{alvarez2020simple,FGLX022,Fre20,Ketch,Kruse}.  In this paper, we focus  on the optimal control of the SIR model \cite{kermack1927contribution} 
\begin{equation}\label{stateeq}
\begin{cases}
\dfrac{ds}{dt}(t)=-\big(\beta-u(t)\big)\, s(t)\, i(t)\\[2ex]
\dfrac{di}{dt}(t)= \big(\beta-u(t)\big)\, s(t)\, i(t)-\gamma i(t)\\[2ex]
s(0)=s_0,\ i(0)=i_0
\end{cases}\qquad ,\,t\in I,
\end{equation}
with $\beta>0$, on an {\em infinite} time horizon $I=(0,\infty)$.
{\color{black} The unknowns $s$ and $i$ in the system \eqref{stateeq} represent the density functions of the susceptible and infectious populations, respectively. The total population is assumed to be constant and the third {\em recovered class} $r$ of SIR is classically obtained by using the conservation of mass. In Section \ref{sec_acs} we recall some other epidemiological concepts related with the SIR model, like, for instance, the {\em herd immunity threshold} and the {\em basic reproduction number} $R_0$.} Under a prescribed threshold $i_M\in(0,1)$, the initial conditions are taken to be $i_0\in(0,i_{M}]$ and $s_0\in(0,1-i_0]$. 
The problem is studied under the control constraints $u\in U:=L^\infty(I;[0,u_{\max}])$ with $0<u_{\max}<\beta$ and the state constraint
$$
i(t)\le i_{M},\ \forall\,t\in I.
$$
Here,  $u_{\max}$ represents the maximal realistic control effort that can be done to avoid the epidemic spread. Also
$i_{M}$ is prescribed and it represents an {\em intensive care unit} (ICU) safety upper bound to the capacity of the health-care system to treat infected patients.

In our problem, $u(t)$ plays the role of a control variable that, typically, models any kind of non-pharmaceutical intervention suitable for reducing the transmission rate $\beta$.    
 The derivatives in
 \eqref{stateeq} are meant in a distibutional sense and the trajectories are constructed from {\em admissible controls} $u\in U$.

The optimal control problem that we aim to consider here is to minimize on the whole time horizon a linear integral cost functional depending both on the state $i$ and on the control $u$, that is 
\begin{equation}\label{cost0}
J^\infty(u,s,i)=\int_I \big(\lambda_1u+\lambda_2i\big)\,dt
\end{equation}
where $(s,i)$ is the epidemic trajectory (i.e.\ solution of \eqref{stateeq}) corresponding to the control $u$. 
The cost functional $J^\infty$ represents the cost of treatments and hospitalization for the populations $i$ of infected individuals 
and its dependence on $u$ allows to capture 
the economic and social cost of non-pharmaceutical interventions like slowdown, isolation, quarantine and distance measures in general. 

For a cost functional depending only on the control (that is for $\lambda_2=0$) 
this problem has been considered in Miclo \cite{Miclo}. When $I$ is a bounded interval, that is for a finite time horizon, it has been studied in \cite{AFG2022}. To be more precise, in \cite{AFG2022} we have considered a finite horizon, state-independent, cost which corresponds to the choice $\lambda_2=0$.  A viability study developed there is useful also in the present paper.  
In particular, we characterized a maximal zone ${\mathcal A}$ in which no control is needed to keep the infections under the level $i_M$,  and a larger zone ${\mathcal B}$ in which policies allowing the trajectory to satisfy the constraints exist. 
Using these characterizations, we established the optimality of the {\em greedy strategy} acting only when the trajectory is about to exit the viable set.  
In particular, it is optimal to stop any control policy as soon as the immunity threshold $s=\frac{\gamma}{\beta}$ is reached. As a consequence, any optimal control for $I=(0,+\infty)$ is optimal also for a finite horizon problem $(0,T)$ with $T$ large enough. 

The greedy policy, which is optimal for the state independent cost considered in \cite{AFG2022},   can no longer be expected to be optimal when the cost \eqref{cost0} also depends on $i$, or,  at least when the contribution of the term involving the state is able to influence the control policy, that is if $\lambda_2$ is big enough.  Moreover, 
differently from the control, the state function $i$ vanishes only as $t\to\infty$ and the reduction to a finite horizon argument used before does not work and must be refined. In particular,  
it will be shown that, in this case as well, a reduction to finite horizon problems can be obtained, albeit this is done by a $\Gamma$-convergence argument.

The variational property of $\Gamma$-convergence implies  that, if $u^T$ is an optimal control (which can be shown to exist) for the problem with finite horizon $T$ (extended by $0$ to the whole infinite horizon) then, as $T$ goes to $\infty$ and up to a subsequence, $u^T$ converges to $u$ weakly* in $L^\infty$,  and $u$ is optimal on the infinite horizon.

Since $u$ is optimal, then it is easy to see that herd immunity must be reached in a finite time  
(see Corollary \ref{corhnc}).  As a consequence, for any $T$ large enough, the susceptible population corresponding to the optimal control $u^T$  satisfies  $s(T)<\frac{\gamma}{\beta}$. We are, then, led to characterize the optimal controls for such $T$ large enough being allowed, in doing it, to use the final condition $s(T)<\frac{\gamma}{\beta}$ as a necessary optimality condition (not as an assumption). 
This property of the approximating optimal controls on a finite horizon has a very important simplifying consequence in the application of Pontryagin's principle, because it implies normality (see the proof of Theorem \ref{ThExtrA} and \ref{B0soc}). For this reason, for this specific problem, considering an infinite horizon is more effective than working on a finite one. Pontryagin's principle, in a version suitable to take into account state contraints, 
allows to characterize the controls $u^T$ and, hence, the control $u$ by taking the weak* limit as $T$ goes to $\infty$.     

The paper is organized as follows. In Section \ref{sec_acs} we state existence, uniqueness and asymptotic behavior of solutions  of the controlled SIR system. In particular, it is shown that if the control has a finite integral on the infinite horizon, then the herd immunity must be reached in finite time.  In Section \ref{EVA} we revisit a viability analysis already done in \cite{AFG2022} by using only elementary methods. This analysis is also useful to enlighten the ``laminar flow argument" that will be used also in the subsequent sections.
In Section \ref{higcs} we introduce the greedy control strategy and show examples of the viability zones compatible with two different times of the Covid-19 epidemic in Italy. In Section \ref{sec_OCP} we formulate the infinite horizon optimal control problem with a general cost and prove the existence of a solution. The case of a linear cost is treated in the subsequent Section \ref{Section3}.  Section \ref{RFHP} is devoted to reduction by $\Gamma$-convergence to a finite horizon problem. In Section \ref{sec:fhphc} we use a version of Pontryagin's theorem, suitable for state constrained problems, to derive necessary conditions of optimality for a finite horizon sufficiently large.          
In Section \ref{SectionExtremalsI} we characterize the structure of the optimal controls on the finite as well as on the infinite horizon.  To numerically illustrate our results we choose the framework of Covid-19 epidemic in Italy and perform some simulations in Section \ref{SectBocop} by using the solver Bocop, \cite{Bocop, BocopExamples}. 

\

\noindent {\bf Notation.} Along the paper, we denote by 
\begin{itemize}
\item $U:=L^\infty(I;K)$ the space of (equivalence classes of) Lebesgue measurable and essentially bounded functions defined on $I$ and taking  values in the compact set $K:=[0,u_{\max}]$;
\item  $Y:=W^{1,\infty}(I;\R^2)$ the Sobolev space of (equivalence classes of) functions that are essentially bounded together with their distributional derivative  defined on $I$ and taking  values in $\R^2$; it is well known that every function in this space has a Lipschitz continuous representative. 
\end{itemize}

\section{Analysis  of the controlled system}\label{sec_acs}

We  refer to the Cauchy problem \eqref{stateeq} as being the {\em controlled system}, because the control function $u$ appears as  an input  in the state differential equations. The uncontrolled one corresponds to the case  in which $u$ is identically equal to $0$.  This section is devoted to study existence, uniqueness and asymptotic behavior of the solutions of the aforementioned controlled system. 
From the point of view of mathematical analysis, the controlled case differs from the uncontrolled one. The uncontrolled system has constant coefficients and this makes it quite easy to study and the solution is continuously differentiable. On the contrary, the controlled system has variable and possibly discontinuous coefficients. Then, $C^1$ solutions are no longer expected and 
their asymptotic behavior is harder to study (see also \cite{Miclo}).   

\begin{theorem}\label{esT}
For every control $u\in U$ and every initial condition
$i_0\in(0,1)$ and $s_0\in(0,1-i_0]$ the controlled system \eqref{stateeq} admits a unique solution $(s,i)\in Y$.
Moreover, 
$(s(t),i(t))\in\mathbb{T}:=\{(s,i)\in(0,1]^2\ :\ s+i\le1\}$ for every $t\in I$.
\end{theorem}

\noindent Let us remark that, in particular, the population densities $s$ and $i$ are always strictly positive and belong to the invariant triangle 
$\mathbb{T}$. A crucial role in the proof is played by the observation that the sum of the two equations gives 
\begin{equation}\label{s+ip}
(s+i)'=-\gamma i.
\end{equation}
Since $i$ is positive, this implies that $s+i$ is decreasing. 

\vspace{2ex}
 
\proof Since the dynamic is locally Lipschitz, then it is classical that we have local existence and uniqueness of an absolutely continuous solution (see for instance \cite[I.3]{Hale1980}). Let us denote by $J$,  with $0\in J$, a time interval in which the unique solution $(s,i)$ exists.  


Let us claim that $s(t)>0$ and $i(t)>0$ for every $t\in J$. Indeed,
considering $i$ as a coefficient, the function $s$ is the unique solution of the linear Cauchy problem
$$
\begin{cases}
s'=-(\beta-u) s i\\
s(0)=s_0>0,
\end{cases}
$$
which is given by
$$
s(t)=s_0\e^{-\int_0^t(\beta-u(\xi)) i(\xi)\,d\xi}\quad\forall\,t\in J.
$$
Hence $s$ is strictly positive in $J$. Similarly, one can prove that also $i$ is strictly  positive. 

\noindent Since, as already remarked, $s+i$ is decreasing, we have
$$
0<s(t)+i(t)\le s_0+i_0\quad\forall\,t\in J,
$$  
and hence $$(s(t),i(t))\in\mathbb{T}:=\{(s,i)\in(0,1]^2\ :\ s+i\le1\}$$ for every $t\in J$.
Being bounded, the solution exists for every $t\in[0,+\infty)$.

Finally, looking at the equations, one sees that $s$ and $i$ have bounded derivatives, hence they are Lipschitz continuous
and the theorem is completely proven. \QED

\begin{remark}[Notation]\label{dsn}{\em
 Whenever the control $u$ and the initial conditions $s(0)=s_0,\ i(0)=i_0$ are fixed,  the unique solution to \eqref{stateeq} will also be denoted by $(s^{s_0,i_0,u}, i^{s_0,i_0,u})$. However, in the sequel, when the initial conditions are fixed or can be easily deduced from the context, we write simply  $(s^{u}, i^{u})$, or even $(s, i)$ if also the control is fixed or easy to deduce. 
This notation will be used  throughout the whole paper. 
 }
 \end{remark}

In addition to the notation just introduced, for the rest of this section, we enforce the assumptions of Theorem \ref{esT}, that is, \[i_0\in(0,1),\ s_0\in(0,1-i_0],\ u\in U.\] 

\begin{lemma}\label{lemec}
The solutions $s$ and $i$ are related by the following 
integro-differential equation
\begin{equation}\label{ibotsub}
s(t)+i(t)-s_0-i_0={\gamma}\int_0^t \frac{s'(\tau)}{(\beta-u(\tau))s(\tau)}\,d\tau
\end{equation}
for every $t\in I$. 
\end{lemma}

\proof 
By integrating \eqref{s+ip} between $0$ and $t$ we have
\begin{equation}\label{ibot}
s(t)+i(t)-s_0-i_0=-\gamma\int_0^t i(\tau)\,d\tau.
\end{equation}
Equation \eqref{ibotsub} is then obtained by using the first state equation to substitute $i$ in the integral. 

\QED

The next theorem deals with the asymptotic behavior of the epidemic trajectory.

\begin{proposition}\label{pcs}
The following propositions hold:
\begin{enumerate}
\item  $s$ is strictly decreasing;
\item $\dis s_\infty:=\lim_{t\to+\infty}s(t)\in(0,1)$ and $\dis i_\infty:=\lim_{t\to+\infty}i(t)=0$;
\item $s_\infty<\frac{\gamma}{\beta-u_{\max}}$.
\end{enumerate}
\end{proposition} 

\proof  {\em 1}. By the first equation and the strict positivity of $s$ and $i$ we have $s'<0$ in $I$; hence the claimed monotonicity. 
 
{\em 2}. By monotonicity and positivity, the limit  $s_\infty:=\lim_{t\to+\infty}s(t)$ exists and is finite. 
Since $s+i$ is decreasing (see \eqref{s+ip}) and positive, then $\lim_{t\to+\infty}\big(s(t)+i(t)\big)$ exists and is finite. 
As a consequence, also the limit
$$
i_\infty:=\lim_{t\to+\infty}i(t)=\lim_{t\to+\infty}\big(s(t)+i(t)-s(t)\big)=\lim_{t\to+\infty}\big(s(t)+i(t)\big)-s_\infty
$$
exists and is finite. 

Let us consider now Lemma \ref{lemec}.
Since  
$\frac{1}{\beta}\le \frac{1}{\beta-u(\tau)}$
and $\frac{s'(\tau)}{s(\tau)}<0$, by using \eqref{ibotsub} we have 
$$
s(t)+i(t)-s_0-i_0\le\frac{\gamma}{\beta}\int_0^t \frac{s'(\tau)}{s(\tau)}\,d\tau=\frac{\gamma}{\beta}\log\frac{s(t)}{s_0}.
$$
Taking the limit as $t\to\infty$ then we get
$$
s_\infty+i_\infty-s_0-i_0\le\frac{\gamma}{\beta}\log\frac{s_\infty}{s_0},
$$
with the convention that $\log0:=-\infty$.
Since  $s$ is positive, we have $s_\infty\ge0$. On the other hand, if we had $s_\infty=0$, by the previous inequality, we would get 
$i_\infty-s_0-i_0\le-\infty$, which is clearly impossible. Hence, $s_\infty>0$.

It remains to prove that $i_\infty=0$. Of course, we have $i_\infty\ge0$. Assume by contradiction that $i_\infty>0$.
Then, there exists $\bar{t}>0$ such that $i(t)>{i_\infty}/{2}$ for all $t\ge\bar{t}$. 
By the first state equation, and since $s(t)\ge s_\infty>0$ and $\beta-u(t)\ge\beta-u_{\max}$, we would have
$$
s'(t)\le -(\beta-u_{\max})s_\infty \frac{i_\infty}{2}<0,\quad\forall\,t\ge\bar{t},
$$
leading to the contradiction $s_\infty=-\infty$. 

{\em 3.} Since $s$ is decreasing, if $s_0\le\gamma/(\beta-u_{\max})$, then there is nothing to prove. 
Let us then assume that $s_0>\frac{\gamma}{\beta-u_{\max}}$. 
Let us denote by
$$
t_{\max}:=\sup\{t\in[0,+\infty)\ :\ s(t)>\dfrac{\gamma}{\beta-u_{\max}}\}
$$
and assume, by contradiction, that $t_{\max}=\infty$. Since $u\le u_{\max}$, we have 
$$
i'(t)=\big((\beta-u(t)) s(t)-\gamma\big)i(t)\ge \big((\beta-u_{\max}) s(t)-\gamma\big)i(t)>0\quad \forall\, t>0,
$$
hence $i$ is increasing. Then, for every $t\in[0,+\infty)$ we have
$$
s'(t)=-(\beta-u(t)) s(t) i(t)\le -(\beta-u_{\max}) s(t) i(t)\le  -\gamma i(t)
\le -\gamma i_0.$$
By integrating on $[0,t)$, we would have  
$$s(t)\le s_0-\gamma i_0 t,\quad \forall\,t\ge0,$$
and therefore $\lim_{t\to+\infty}s(t)=-\infty$, which contradicts 2.

It follows that  $t_{\max}<\infty$. Then, by continuity, $s(t_{\max})=\frac{\gamma}{\beta-u_{\max}}$ and the claim follows by the strict
monotonicity of $s$.  

\QED

By the second equation of the uncontrolled SIR system \eqref{stateeq} (with $u=0$), since $i$ is positive, we deduce
that 
$$
i'(t)\le 0\ \iff\ s(t)\le\frac{\gamma}{\beta}.
$$
The number $\frac{\gamma}{\beta}$ is the {\em herd immunity threshold}. Since $s$ is decreasing,
when immunity is reached, $i$ starts to decrease and continues to do it for every time after. {\color{black} In the SIR model, the immunity threshold is the reciprocal of the {\em basic reproduction number $R_0$} (see, for instance, \cite{Heth}).}

\begin{remark}\label{siuc}{\em
In Proposition \ref{pcs} we have proved that $s_\infty<\frac{\gamma}{\beta-u_{\max}}$. 
This number is strictly greater than $\frac{\gamma}{\beta}$. If $u\equiv 0$, then we have $s_\infty<\frac{\gamma}{\beta}$, by simply taking $u_{\max}=0$ in point {\em 3.\ }of Proposition \ref{pcs}.
On the other hand, it is easy to see that  the herd immunity may not be reached 
if the control $u$ is big enough. That is, it may happen that $s_\infty>\frac{\gamma}{\beta}$.}
\end{remark}

The next theorem provides a very useful sufficient condition that ensures that herd immunity is reached in finite time.
   Proposition \ref{pcs} is quite standard, even if not trivial for the controlled system. On the contrary, to the best of our knowledge, Theorem \ref{uL1} is new. It plays an important role in the optimal control problems that will be considered in the sequel. 

\begin{theorem}\label{uL1}
If  $u\in L^1(I)$ then $s_\infty<\frac{\gamma}{\beta}$.
\end{theorem}

The following lemma will be used in the proof. 

\begin{lemma}\label{notn}
There exists $n_0\in\N$ such that 
$$
\forall\,n\ge n_0\ \exists\,t_n>0,\, t_n\to\infty\ :\ i(t_n)=\frac{1}{n}\mbox{ and }i(t)\le\frac{1}{n}\ \forall\,t\ge t_n.
$$
\end{lemma}

\proof Let $n_0>\frac{1}{i_0}$. Let us define
$$
t_{n_0}:=\sup\{t\ge0\ :\ i(t)=\frac{1}{n_0}\}.
$$ 
Since $i_\infty=0$ and $i$ is continuous, it follows that $t_{n_0}<\infty$. Moreover, $i(t_{n_0})=\frac{1}{n_0}$  and 
 $i(t)\le \frac{1}{n_0}$ for every $t\ge t_{n_0}$.

By continuity, there exists $t> t_{n_0}$ such that $i(t)=\frac{1}{n_0+1}$. Let us define
$$
t_{n_0+1}:=\sup\{t> t_{n_0}\ :\ i(t)=\frac{1}{n_0+1}\}.
$$ 
As before, we have $t_{n_0+1}<\infty$, $i(t_{n_0+1})=\frac{1}{n_0+1}$ and 
 $i(t)\le \frac{1}{n_0+1}$ for every $t\ge t_{n_0+1}$.

We construct the rest of the sequence by induction on $k\in \N$ as follows. Assume that $t_{n_0+k}$ be such that  
$t_{n_0+k}<\infty$, $i(t_{n_0+k})=\frac{1}{n_0+k}$ and 
 $i(t)\le \frac{1}{n_0+k}$ for every $t\ge t_{n_0+k}$. 
By continuity, there exists $t> t_{n_0+k}$ such that $i(t)=\frac{1}{n_0+k+1}$. Let us define
$$
t_{n_0+k+1}:=\sup\{t> t_{n_0+k+1}\ :\ i(t)=\frac{1}{n_0+k+1}\}.
$$ 
As before, we have $t_{n_0+k+1}<\infty$, $i(t_{n_0+k+1})=\frac{1}{n_0+k+1}$ and 
 $i(t)\le \frac{1}{n_0+k+1}$ for every $t\ge t_{n_0+k+1}$.

By construction, the sequence $(t_n)_{n\ge n_0}$ is increasing. It remains only to prove that $t_n\to\infty$. Assuming by contraction that, instead, $t\to t_\infty<\infty$ then we would have $i=0$ for every $t>t_\infty$ (indeed, for every $t>t_\infty$ we would have $i(t)<1/n$ for every $n\ge n_0$). But this is impossible since $i(t)>0$ for every $t\in[0,+\infty)$. The lemma is completely proved.  
\QED

\proof of Theorem \ref{uL1}. Let us assume, by contradiction, that $s_\infty\ge\frac{\gamma}{\beta}$.

Let $n_0\in\N$ and $(t_n)_{n\ge n_0}$ be as in Lemma \ref{notn}.
Since $i_\infty=0$, we have
\begin{equation}\label{itnui}
-i({t_n})=\int_{t_n}^\infty i'\,dt=\int_{t_n}^\infty[(\beta-u)s-\gamma]i\,dt\ge-\frac{\gamma}{\beta}\int_{t_n}^\infty ui\,dt
\end{equation}
where, in the last inequality, we used that fact that, since $s$ is decreasing and $s_\infty\ge\frac{\gamma}{\beta}$, then  $s(t)>\frac{\gamma}{\beta}$ for every $t\in[0,+\infty)$. On the other hand, since $i(t)\le\frac{1}{n}$ for every $t\ge t_n$, we get
\begin{equation*}
-\frac{\gamma}{\beta}\int_{t_n}^\infty ui\,dt\ge -\frac{\gamma}{\beta n}\int_{t_n}^\infty u\,dt.
\end{equation*}
By putting this inequality together with \eqref{itnui} and using the fact that $ni({t_n})=1$, then we get
$$
\int_{t_n}^\infty u\,dt\ge \frac{\beta}{\gamma},\quad \forall\,n\ge n_0.
$$
Since $t_n\to\infty$, this implies 
$$
\int_{0}^\infty u\,dt=\infty,
$$
which contradicts the assumption $u\in L^1(I)$.

 \QED

\begin{remark}\label{ufs}{\em The map $s:[0,+\infty)\to(s_\infty,s_0]$ 
is strictly decreasing, hence invertible and Lipschitz continuous together with its inverse $t=t(s)$.
Then it can be used to perform a change of variable in the Lebesgue integral in \eqref{ibotsub} (see for instance 
\cite{H93}) and we get
\begin{equation}\label{ibotsubs}
s+i(s)-s_0-i_0={\gamma}\int_{s_0}^s \frac{d\sigma}{(\beta-u(\sigma))\sigma}
\end{equation}
for every $s\in[s_\infty,s_0]$. Let us note that here, with a small abuse of notation, $u$ and $i$ are considered as functions of $s$, through the composition with the function $t(s)$ (that is, $u(\sigma)$ stays for $u(t(\sigma))$ and the same holds for $i$). In this sense, the control $u$ is defined only in the interval $(s_\infty,s_0]$, 
and equality \eqref{ibotsubs} does not make sense for $s\in[0,s_\infty)$.

When the control $u$ is constant, the integral at the right hand side of \eqref{ibotsubs} (or \eqref{ibotsub}) can be explicitly computed, leading to the following relation 
\begin{equation}\label{iscc}
i=s_0+i_0-s+\frac{\gamma}{\beta-u}\big(\log s-\log{s_0}\big)
\end{equation}
for every constant control $u\in[0,u_{\max}]$. Differently from the previous one, this equality makes sense for every $s\in(0,1]$ because the constant control $u$ can be thought to be defined everywhere. Of course, when $s$ is taken in $(0,s_\infty)$, 
it returns a negative value of $i$, which is not physically relevant. Furthermore, the equality \eqref{ibotsubs} itself can be extended to $s\in(0,1]$ by extending the control, if needed. 
 }
\end{remark}

\section{An elementary viability analysis}\label{EVA}

The solution of the controlled system 
$$
\begin{array}{rcl}
I&\to&\R^2\\
t&\mapsto&\big(s(t),i(t)\big)
\end{array}
$$
is a curve (the {\em epidemic trajectory}) in the plane $si$.  Since $s$ is a strictly decreasing function of time, the epidemic trajectory is graph of a function $i=i(s)$; the curve is  traveled in the direction of decreasing $s$. 

\begin{figure}[H]\label{etcc_fig}
\begin{center}
\includegraphics[width=0.7\textwidth]{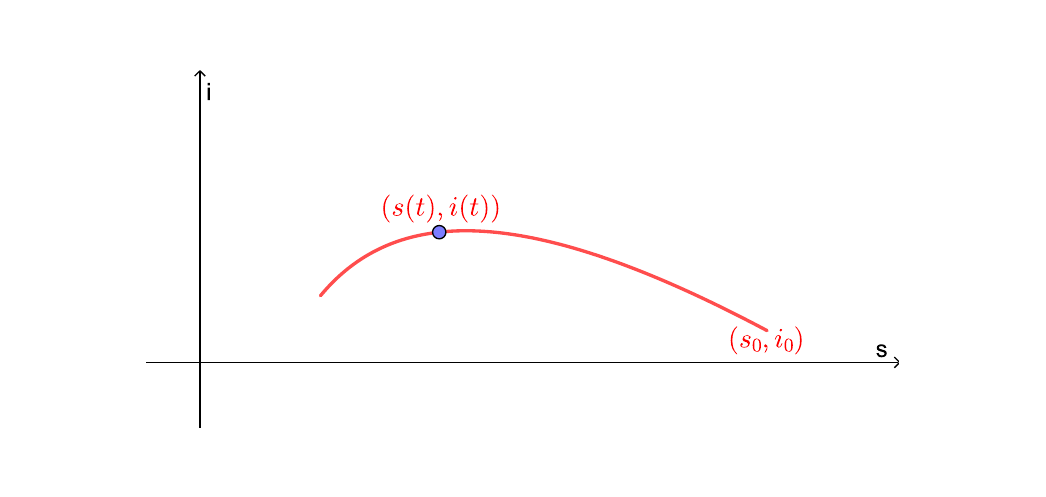} 
\end{center}
\caption{The epidemic trajectory}
\end{figure}

The following viability zones have been already defined in \cite{AFG2022}.

\begin{definition}
The set 
$$
\mathcal{B}:=\{(s_0,i_0)\in\mathbb{T}\ :\exists\,u\in U\mbox{ s.t.\ }i^{s_0,i_0,u}(t)\leq i_M,\ \forall\,t\in I\}
$$
is called {\em feasible} or {\em viable} set. 

The set 
$$
\mathcal{A}:=\{(s_0,i_0)\in\mathbb{T}\ : i^{s_0,i_0,0}(t)\leq i_M,\ \forall\,t\in I\}
$$
is called {\em no-effort} or {\em safe} zone. 
\end{definition}

\noindent The viable set $\mathcal{B}$ is the maximal set of initial configurations on which at least one trajectory satisfies the constraint. When the epidemic initial state $(s_0,i_0)$ is inside $\mathcal{B}$, the epidemic is ``under control"'. The safe zone $\mathcal{A}$ is the maximal set of initial configurations such that the associated uncontrolled trajectories satisfy the constraint. These sets play an important role in the synthesis of optimal control strategies.

\

\noindent The epidemic trajectories corresponding to constant controls are graphics $i=i(s)$ given by \eqref{iscc}. Let us note that these 
functions are defined for every $s>0$, while the epidemic trajectory is only the part of the graph corresponding to $s\in(s_\infty,s_0]$.  

They are easy to study, and play a crucial role in the viability analysis. Given a prescribed (fixed) constant control $u$, 
and as $(s_0,i_0)$ varies, the family of graphs $\{s\mapsto i^{s_0,i_0,u}(s)\}_{(s_0,i_0)}$
satisfies the following properties (very easy to check).
\begin{enumerate}
\item Any graph of the family is strictly concave.
\item The graphs have all a unique maximum point in $\frac{\gamma}{\beta-u}$ (independent of $s_0$ and $i_0$).
\item Any graph goes from $(s_0,i_0)$ to $(s_\infty,0)$ in infinite time (because $i$ is always strictly positive and $i_\infty=0$).
\item It is a {\em laminar flow}, that is, the graphs of the family are, two by two, disjoint or identically equal (recall that the control is constant and fixed). See the figure below for two laminar flows with different controls. 
\end{enumerate}
The figure shows some epidemic trajectories with constant controls. Of course only the arcs contained in the half plane of positive $i$ 
have physical meaning.

\begin{figure}[H]\label{et0max_fig}
\begin{center}
\includegraphics[width=0.8\textwidth]{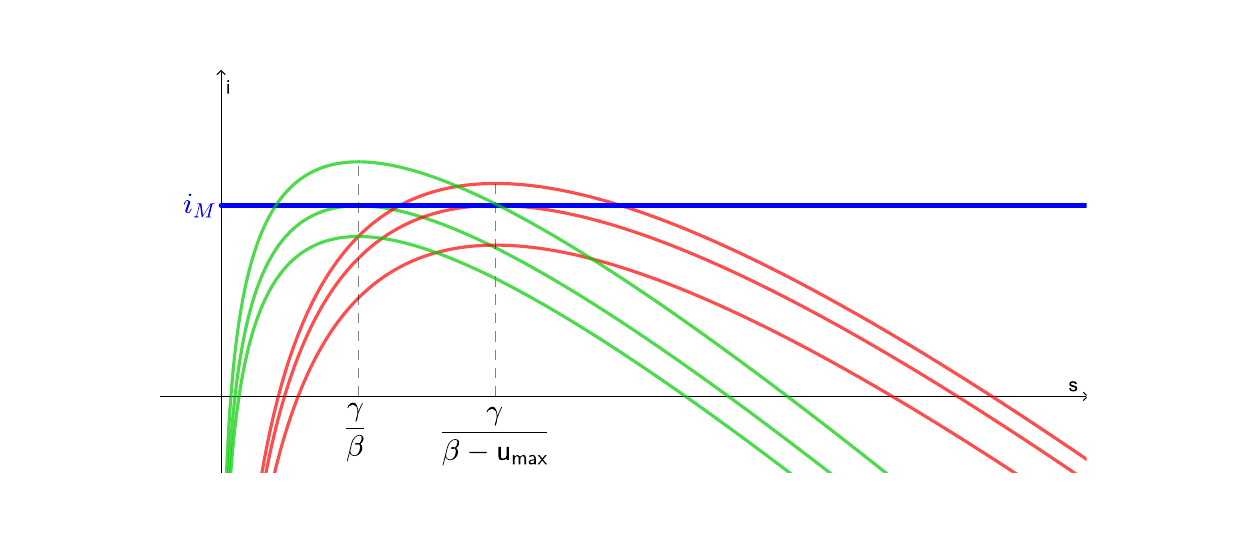}\\
\end{center}
\caption{Epidemic trajectories with $u=0$ (in green) and $u=u_{\max}$ (in red)}
\end{figure}
\noindent The picture above suggests how to construct two relevant curves, $\Phi_{0}$ and $\Phi_{max}$, in black in the figure below.
\begin{figure}[H]\label{phimax_fig}
\begin{center}
\includegraphics[width=0.8\textwidth]{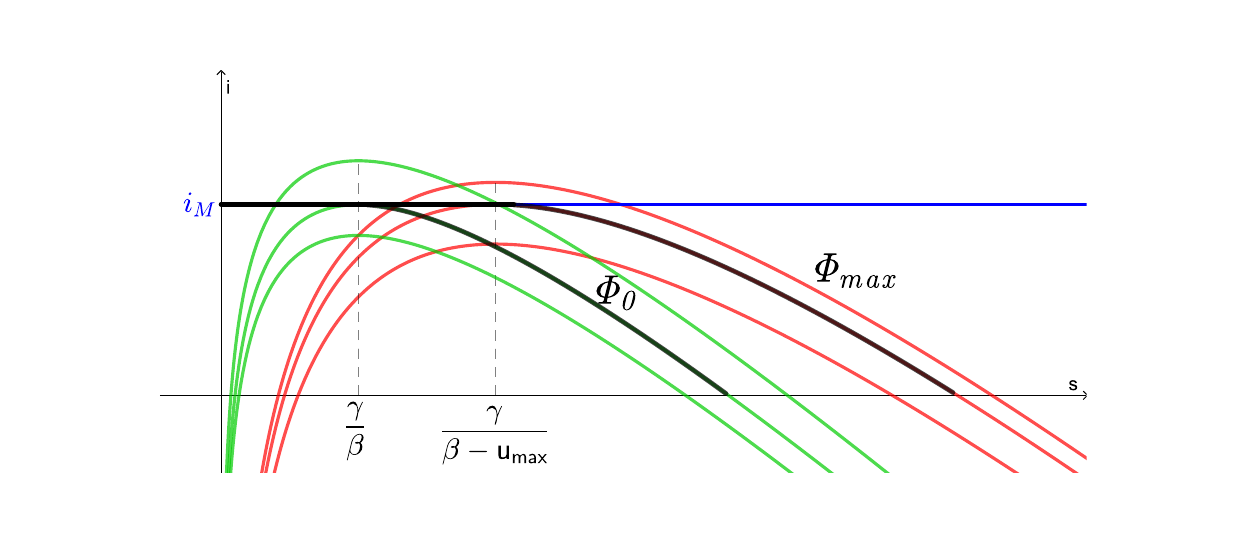}
\end{center}
\caption{The curves $\Phi_{0}$ and $\Phi_{\max}$}
\end{figure}
\noindent They are obtained by glueing together the constant $i=i_M$ with the epidemic trajectories with $u=0$ and $u=u_{\max}$ and passing trough the maximum points
$$\big(\frac{\gamma}{\beta},i_M\big),\qquad \big(\frac{\gamma}{\beta-u_{\max}},i_M),$$
respectively, 
that is
\begin{eqnarray*}
\Phi_{0}(s)&=&\begin{cases}
i_M&\mbox{if } 0<s<\frac{\gamma}{\beta},\\
\frac{\gamma}{\beta}+i_M-s+\frac{\gamma}{\beta}\big(\log s-\log{\frac{\gamma}{\beta}}\big)&\mbox{if } s\ge \frac{\gamma}{\beta},
\end{cases}\\
\Phi_{\max}(s)&=&\begin{cases}
i_M&\mbox{if } 0<s<\frac{\gamma}{\beta-u_{\max}},\\
\frac{\gamma}{\beta-u_{\max}}+i_M-s+\frac{\gamma}{\beta-u_{\max}}\big(\log s-\log{\frac{\gamma}{\beta-u_{\max}}}\big)&\mbox{if } s\ge \frac{\gamma}{\beta-u_{\max}}.
\end{cases}
\end{eqnarray*}
The region under the curve $\Phi_{max}$ is a viable set. 
Indeed, by choosing a point inside it and applying the maximum control $u=u_{\max}$ the corresponding trajectory
stays always on, or strictly under, the curve $\Phi_{max}$ (due to the laminar flow property 4.; see the red curve under $\Phi_{\max}$ in Figure \ref{phimax_fig}).  
If, instead, we start outside this region then the maximum control is not enough to take the epidemic curve under the 
level $i_M$ (always due to property 4.). 
{\color{black} To complete a formal proof of the claimed viability, it remains to prove that this last fact holds also for any other control (not only for $u=u_{\max}$).} That is, the following characterization of the viable set holds.

\begin{theorem}
$\mathcal{B}=\{(s,i)\in\mathbb{T}\ :\ 0<i\le\Phi_{\max}(s)\}.$
\end{theorem}

\proof As we already remarked, the inclusion
$$
\{(s,i)\in\mathbb{T}\ :\ 0<i\le\Phi_{\max}(s)\}\subseteq\mathcal{B}
$$
is trivially proven  by choosing $u=u_{\max}$, and it remains only to prove the opposite inclusion.

Let us then consider $(s_0,i_0)\in\mathbb{T}$ and $u\in U$ be such such that $i^{s_0,i_0,u}\leq i_M$ for all $t\in I$.
We have only to consider the case in which 
$s_0>\frac{\gamma}{\beta-u_{\max}}$ since, otherwise,  $(s_0,i_0)\in\{(s,i)\in\mathbb{T}\ :\ 0<i\le\Phi_{\max}(s)\}$ and the inclusion is proved.
We aim to prove that $i_0\le \Phi_{\max}(s_0)$, that is
\begin{equation}\label{tfm}
i_0\le \frac{\gamma}{\beta-u_{\max}}+i_M-s_0+\frac{\gamma}{\beta-u_{\max}}\big(\log s_0-\log{\frac{\gamma}{\beta-u_{\max}}}\big).
\end{equation}
By Lemma \ref{lemec}, denoting by $i=i^{s_0,i_0,u}$ and $s=s^{s_0,i_0,u}$, we have
$$
i_0=i(t)-s_0+s(t)-{\gamma}\int_0^t \frac{s'(\tau)}{(\beta-u(\tau))s(\tau)}\,d\tau,\ \forall t\ge0,
$$
for every $t\ge0$. Since $i\le i_M$, 
by using $\frac{1}{\beta-u(\tau)}\le \frac{1}{\beta-u_{\max}}$
and $\frac{s'(\tau)}{s(\tau)}<0$ to estimate the integral, we get
$$
i_0\le i_M-s_0+s(t)-\frac{\gamma}{\beta-u_{\max}}\big(\log s(t)-\log s_0\big). 
$$
which holds for every $t\ge0$. Since, by point 3. of Proposition \ref{pcs} and by continuity, we have that there
exists $t_{\max}\ge0$ such that  $s(t_{\max})=\frac{\gamma}{\beta-u_{\max}}$, then the evaluation of the previuos inequality in $t=t_{\max}$ gives
$$
i_0\le i_M-s_0+\frac{\gamma}{\beta-u_{\max}}-\frac{\gamma}{\beta-u_{\max}}\big(\log\frac{\gamma}{\beta-u_{\max}}-\log s_0\big)
$$
which is exactly \eqref{tfm}. \QED

\vspace{2ex}

Similarly, the region under the curve $\Phi_{0}$ is a safe zone. 
Indeed, by choosing a point inside it and applying the control $u=0$ the corresponding trajectory
stays always on, or strictly under, the curve $\Phi_0$ (see the green curve under $\Phi_{0}$ in Figure \ref{phimax_fig}).  
If, instead, we start outside this region then the control $u=0$ is not enough to take the epidemic curve under the 
level $i_M$,
that is, the following characterization of the safe zone holds.

\begin{theorem}
$\mathcal{A}=\{(s,i)\in\mathbb{T}\ :\ 0<i\le\Phi_{0}(s)\}.$
\end{theorem}

%

\begin{remark}{\em
The previous theorems have been proven in \cite[Theorem 2.3]{AFG2022} by using viability tools. The proofs given here are elementary.
}
\end{remark}

\section{Greedy control strategy and examples}\label{higcs}

Starting from a point inside the viability set $\mathcal{B}$ and parametrizing by $s$ instead of $t$ (see Remark \ref{ufs}), the {\em greedy control strategy} consists in using for every $s$ the minimal control effort that takes the system inside the set  $\mathcal{B}$; namely,
\begin{enumerate}
\item 
$u_g=0$  till the associated trajectory reaches a point on the curve $\Phi_{\max}$ 
and, afterwards, follow this curve till reaching herd immunity by using the following strategy
\item $u_g=u_{\max}$, as long as $s>\frac{\gamma}{\beta-u_{\max}}$,
\item $u_g=\beta-\frac{\gamma}{s}$,  as long as $s\in(\frac{\gamma}{\beta},\frac{\gamma}{\beta-u_{\max}})$,
\item $u_g=0$,  for every $s\ge\frac{\gamma}{\beta}$, that is, once that herd immunity has been reached.
\end{enumerate}  
\begin{figure}[H]
\begin{center}
\includegraphics[width=0.85\textwidth]{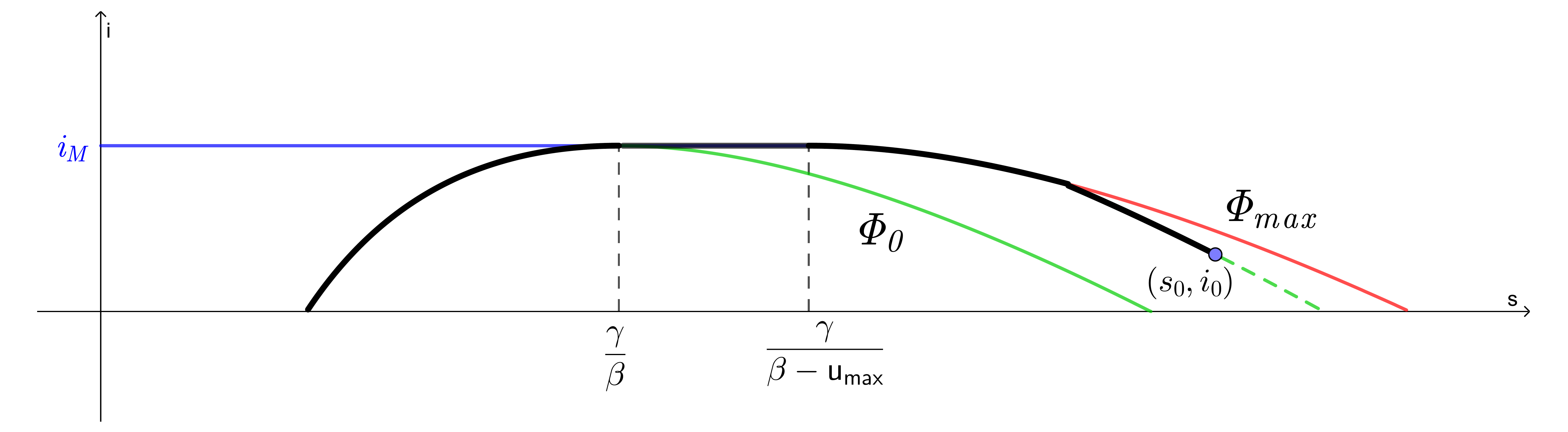}
\end{center}
\caption{The greedy strategy (in black) starting from $(s_0,i_0)\in\mathcal{B}\setminus\mathcal{A}$}
\label{gs}
\end{figure}
In particular, by applying the control 3.\ we have $i'=0$ (by the second state equation), hence $i$ is taken at the level $i_M$ till to reach the herd immunity.
Note that such control is admissible because it is an increasing function of $s$ on the interval $(\frac{\gamma}{\beta},\frac{\gamma}{\beta-u_{\max}})$ taking its minimum value $0$ in the left point $\frac{\gamma}{\beta}$ and maximum value  $u_{\max}$ in the right point $\frac{\gamma}{\beta-u_{\max}}$.

\begin{figure}[H]
\begin{center}
\includegraphics[width=0.9\textwidth]{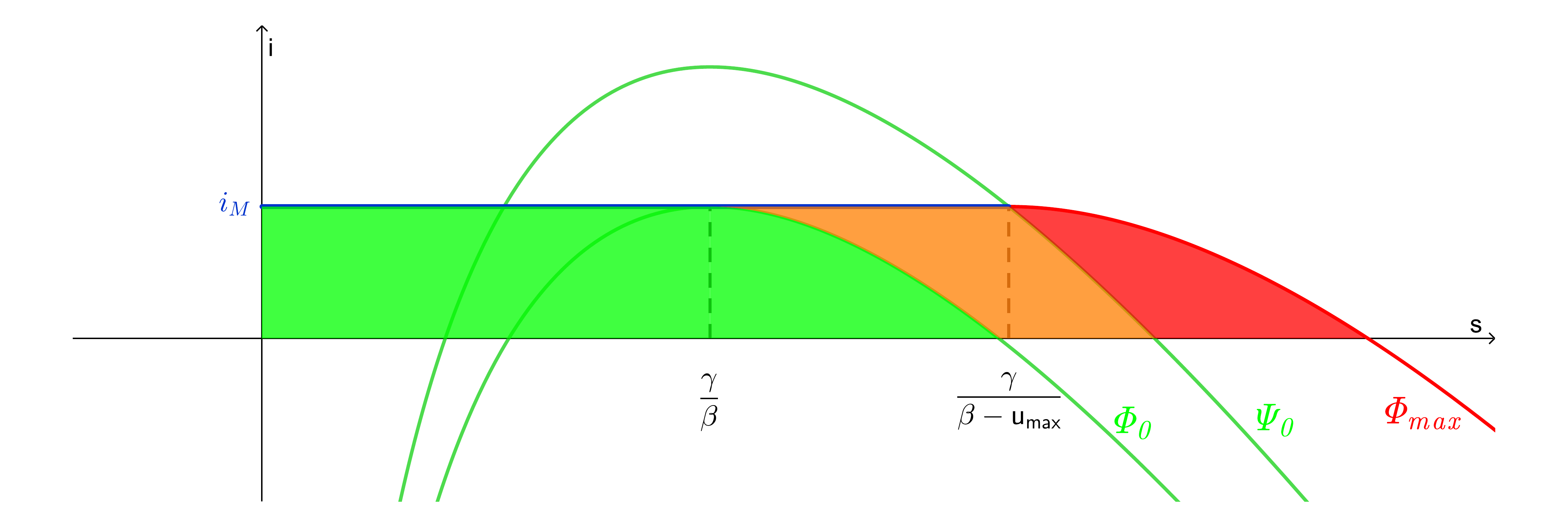}
\end{center}
\caption{The curve $\Psi_0$}
\label{psi0_fig}
\end{figure}

\noindent A third important curve is the graph $\Psi_0$ displayed in Figure \ref{psi0_fig}.

\noindent It is a curve travelled by control $u=0$ and passing trough the point $(\frac{\gamma}{\beta-u_{\max}},i_M)$ and truncated at the level $i_M$, that is
$$
\Psi_{0}(s)=\begin{cases}
i_M&\mbox{if } 0<s<\frac{\gamma}{{\beta-u_{\max}}},\\
\frac{\gamma}{\beta-u_{\max}}+i_M-s+\frac{\gamma}{\beta}\big(\log s-\log{\frac{\gamma}{\beta-u_{\max}}}\big)&\mbox{if } s\ge \frac{\gamma}{{\beta-u_{\max}}}.
\end{cases}
$$
It has the property that, starting from an initial point $(s_0,i_0)$ between $\Phi_0$ and $\Psi_0$ with control $u=0$,
the level $i_M$ is reached  without using the maximum control $u_{\max}$, that is, without imposing a lockdown.
In the picture below, the safe zone $\mathcal{A}$ is coloured in green. The red zone corresponds to the points of $\mathcal{B}$ laying between the curves $\Psi_0$ and $\Phi_{\max}$ where a lockdown (i.e.\ $u=u_{\max}$) is needed.

\begin{remark}{\em
The greedy strategy has been proved to be optimal (in \cite[Theorem 5.6]{AFG2022}) for the case of a cost functional depending only on the control (that is $\lambda_2=0$).}
\end{remark}

\begin{remark}\label{hirft}{\em
Along the greedy strategy, the herd immunity is always reached in finite time.  
This can 
be easily seen by observing that, starting from $(s_0,i_0)\in \mathcal{B}\setminus\mathcal{A}$, the curve $\Phi_{\max}$ is reached along a curve with control $u=0$ in finite time (indeed, the time to reach would be infinite only if the intersection point was $(s_\infty^{s_0,i_0,0},0)$, that is on the $s$ axis). Analogously, the time spent along the red part of the curve $\Phi_{\max}$ is finite or null. Finally, also the time eventually spent along the line $i=i_M$ can be computed (as well as in the other cases, by the way) by using the fact that, with the control given by 2.\ the first SIR equation becomes
$$
\frac{ds}{dt}=-\gamma i_M,
$$
and, hence, the time spent along the line is, at most,
$$
T_{\rm stab}=\int_{\frac{\gamma}{\beta-u_{\max}}}^\frac{\gamma}{\beta}\frac{dt}{ds}\,ds= 
-\int_{\frac{\gamma}{\beta-u_{\max}}}^\frac{\gamma}{\beta}\frac{1}{\gamma i_M}ds=
\frac{u_{\max}}{\beta(\beta-u_{\max})i_M}.
$$ 
The other times can be computed in an analogous way. 
   }
\end{remark}

 \begin{example}\label{ex_covita}{\em 
 To show examples we have to choose the epidemic parameters $\beta$, $\gamma$, $u_{\max}$, 
 and the ICU threshold $i_M$. In this example we do it according to the available Covid-19 data in Italy at different times of the epidemic. In particular, we focus on the initial stage of March 2020 (mainly in Northern Italy) and on autumn 2021 when the government updated the control policies against the diffusion of the Delta variant.

According to  \cite{gea021} (see also \cite{Talice068302})  which reports a range between $2.83$ and $3.10$, we take the basic reproduction number $R_0=3$. This value is a little bit more than the arithmetic average and takes into account the fact that the first part of the epidemic involved mainly the northern's regions.

Moreover, by considering a time-to-recovery of 14 days (\cite{BMDetal020}), we get $\gamma=\frac{1}{14}\simeq 0.0714$ and hence $\beta=R_0\gamma\simeq 0.2142$.

 About the numbers ICU places, $p_{\rm icu}$, an optimistic estimate  is of $14$ over $100\,000$ inhabitants, which is, in fact, the Italian target after Decreto-Legge 19.05.2020 \cite{DL19052020}. The reported critical cases in March 2020 was $c_{\rm icu}=4.5\%$ (see \cite{ISS_ig25march2020}). Hence,  $i_M=\frac{p_{\rm icu}}{c_{\rm icu}}=0.0031$. 
 The maximal control effort that we assume to be done is
 $u_{\max}=0.63\beta$, which corresponds to $63\%$ reduction of the transmission rate, according to the impact estimate of the lockdown on trasmissions given in  \cite{gea021}.

 The figure shows 
 the viability zones for the previuos choice of the epidemic parameters (summarizing: $\gamma=0.0714$, $\beta=0.2142$, $i_M=0.0031$, $u_{\max}=0.135$) with an initial point $(s_0,i_0)=(0.94,0.001)$ that will be considered in the subsequent Example \ref{ita_cov_noICU}.

 \begin{figure}[H]
\begin{center}
\includegraphics[width=0.8\textwidth]{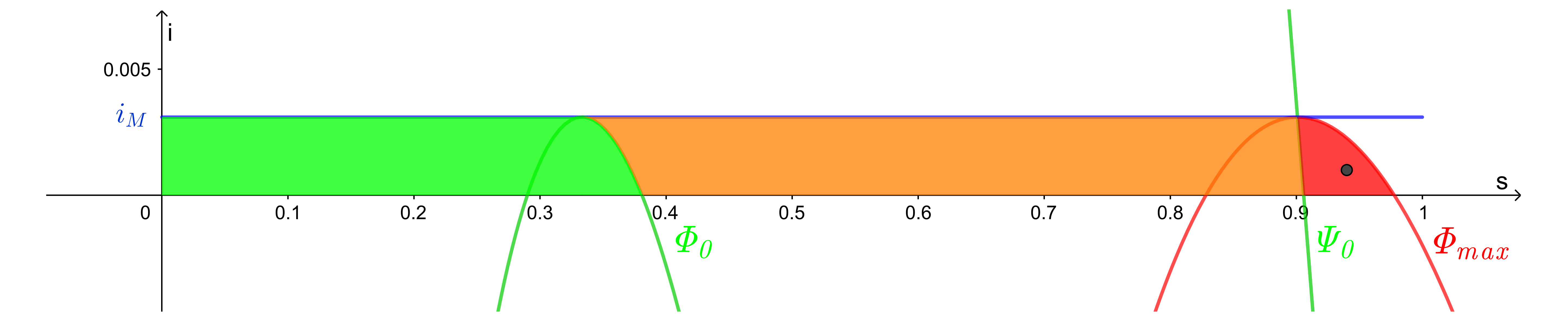}
\end{center}
\caption{The viability zones at the initial stage of Covid-19 epidemic in Italy}
\label{ita_ini}
\end{figure}

\noindent In 2021 the Delta variant was more than two times transmissible than the ancestral Sars-CoV-2 virus (see for instance \cite{taab124}). 
 Accordingly, in this case, we choose $\beta=0.5$. With DL 23 Luglio 2021, n. 105, \cite{DL23072021}, the Italian government setted the  alert level at 150 cases per $100\,000$ inhabitants. Since the time to recovery accredited by WHO was always (in mean) of
$14$ days, we can estimate that the health system was considered able to support a number of infections $14$ times bigger, that is  
$$i_{M}=14\frac{150}{100\,000}=0.021$$
possibly due also to the effect of the vaccinations that reduced the fraction of critical cases.

By updating data in this way (that is, by taking: $\gamma=0.0714$, $\beta=0.5$, $i_M=0.021$, $u_{\max}=0.315$),
we obtain the following figure.

 \begin{figure}[H]
\begin{center}
\includegraphics[width=0.95\textwidth]{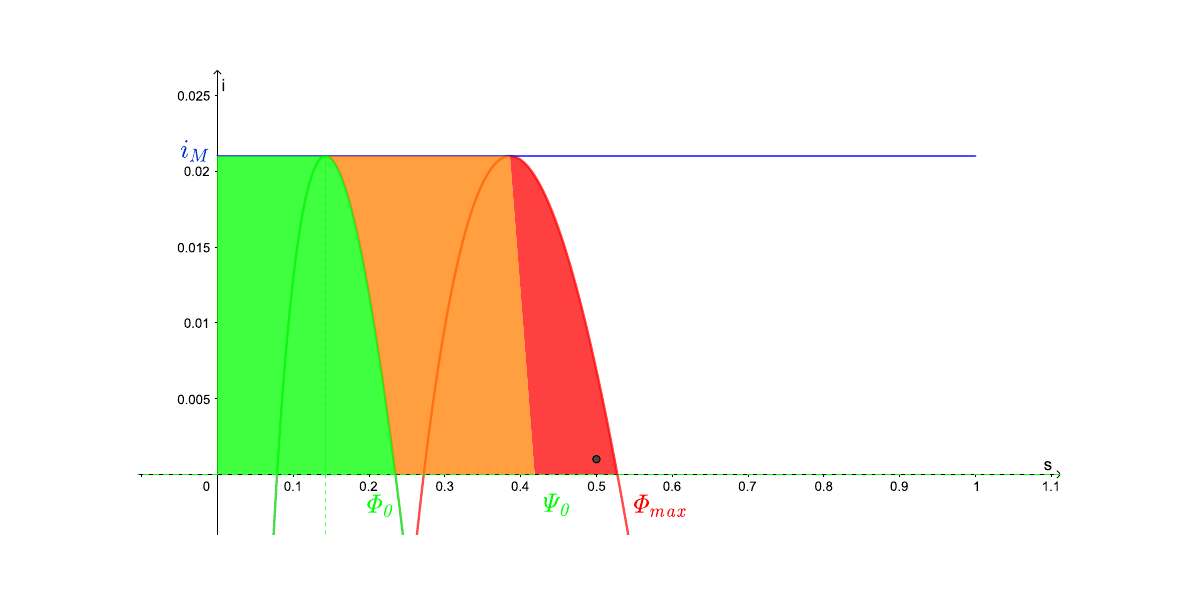}
\end{center}
\caption{Viability zones for the Delta variant (autumn 2021) in Italy}
\label{ita_aut}
\end{figure}
}
\end{example}

\noindent The figure shows also the initial point $(s_0,i_0)=(0.5,0.001)$ that will be considered in the subsequent Example \ref{it2}.

\section{The optimal control problem for a general cost}\label{sec_OCP}

The optimal control problem consists in minimizing a {\em cost functional} $J:  U\times Y\to[0,\infty]$
of the form 
\begin{equation}\label{cost}
J(u,s,i)=\int_0^{\infty}f_0(t,s,i,u)\,dt,
\end{equation}
where $f_0$ is a given {\em running cost} such that the integral makes sense, over the
set of {\em state equations} \eqref{stateeq} with initial conditions $(s_0,i_0)$ in the viable set $\mathcal{B}$ and 
under the {\em ICU constraint} on the trajectories of \eqref{stateeq}
\begin{equation}\label{icug}
i(t)\le i_M,\quad \forall t\in[0,+\infty).
\end{equation}
An {\em optimal solution} to the control problem \eqref{cost}-\eqref{icug}-\eqref{stateeq} is a vector function 
$(u,s,i)\in L^\infty(I;K)\times W^{1,\infty}(I;\R^2)$ that minimizes the cost $J$ and satisfies the set of state equations and the upper bound on $i$. The function $u$ is an {\em optimal control} and $(s,i)$ an {\em optimal state or trajectory}. Accordingly, $U=L^\infty(I;K)$ is the {\em space of controls} and $Y=W^{1,\infty}(I;\R^2)$ is the {\em space of states}.

The following semicontinuity and existence theorem for a very general cost functional holds.

\begin{theorem}\label{gen_ex_th}
If $f_0:I\times\R^2\times\R\to[0,\infty)$ is a normal convex integrand, that is 
it is measurable with respect to the Lebesgue $\sigma$-algebra on $I$ and the Borel $\sigma$-algebra on $\R^2\times\R$ 
and there exists a subset $N\subset I$ of Lebesgue measure zero such that  
\begin{enumerate}
\item $f_0(t,\cdot,\cdot,\cdot)$ is lower semicontinuous for every $t\in I\setminus
N$,
\item $f_0(t,s,i,\cdot)$ is convex for every $t\in I\setminus
N$ and $s,i\in\R$,
\end{enumerate}
then 
\begin{enumerate}
\item[a.] the cost functional $J$ defined by \eqref{cost} is weakly* lower semicontinuous,
\item[b.] 
there exists an optimal solution $(u,s,i)$ to the control problem \eqref{cost}-\eqref{stateeq}.
\end{enumerate}
\end{theorem}

\noindent To prove the existence of an optimal solution we  observe that it is equivalent to prove the existence of a minimizer
of the functional 
\begin{equation}\label{minF}
F(u,s,i):=J(u,s,i)+\chi_{\Lambda}(u,s,i)+\chi_{i\le i_M}(i)
\end{equation}
where $\Lambda$ is the set of admissible pairs, that is all state-control vectors $(u,s,i)$ that satisfy the initial value problem \eqref{stateeq}, while $\chi_\Lambda $ denotes the indicator function of $\Lambda$ that takes the value $0$ on $\Lambda$ and $+\infty$ otherwise; similarly, the function $\chi_{i\le i_M}(i)$ is $0$ if $i(t)\le i_M$ for every $t\in I$, and $+\infty$ otherwise. 

\vspace{1ex}
 
\begin{proof} By De Giorgi and Ioffe's Theorem (see for instance \cite[Theorem 7.5]{FL_MMCVLp} or \cite[Section 2.3]{Buttazzo89}), for 
every $T\in(0,+\infty)$ the cost functional 
$J_T:  L^\infty((0,T);K)\times W^{1,\infty}((0,T);\R^2)\to[0,\infty]$ 
defined by
\begin{equation*}
J_T(u,s,i)=\int_0^{T}f_0(t,s,i,u)\,dt
\end{equation*}
is  weakly* lower semicontinuous. Then we obtain easily that also $J$ is lower semicontinuous. Indeed,
given weakly* converging sequences $(s_n,i_n,u_n)\weaks(s,i,u)$, since the integrand is non-negative, we have
\begin{eqnarray*}
\int_0^{\infty}f_0(t,s,i,u)\,dt&=&\sup_{T>0}\int_0^{T}f_0(t,s,i,u)\,dt\\
&\le&
\sup_{T>0}\liminf_{n\to\infty}\int_0^{T}f_0(t,s_n,i_n,u_n)\,dt\\
&\le& \liminf_{n\to\infty}\int_0^{\infty}f_0(t,s_n,i_n,u_n)\,dt.
\end{eqnarray*}
This proves {\em a.} Let us now prove {\em b.} 
On the domain of $F$, that is the space $L^\infty(I;K)\times W^{1,\infty}(I;\R^2)$, we consider the topology given by the product of the weak* topologies of the two spaces and aim to prove sequential lower semicontinuity and coercivity  
of the functional $F$ with respect to this topology. By the 
Direct Method of the Calculus of Variations {\color{black}(see, for instance, Buttazzo \cite[Sec.\ 1.2]{Buttazzo89})}, these properties imply the existence of a solution to the minimum problem. They are direct consequences of the  
fact that the space of controls is weakly* compact, that the assumptions on $f_0$ imply that the cost functional $J$ is weakly* lower semicontinuous 
 and the fact that the sets $\Lambda$ {\color{black} and $\{i\le i_M\}$ are} closed with respect to the weak* convergence. 
{\color{black} The claimed closedness of such sets follows by the application of Rellich compactness theorem, which ensures that weakly* converging sequences in $W^{1,\infty}(I)$ are, up to subsequences, uniformly converging on every bounded subset $J\subset I$ (see for instance \cite[Theorem 8.8 and Remark 10]{Brezis2011}).}  

\QED
\end{proof}

\begin{remark}{\em 
The requirement on $f_0=f_0(t,s,i,u)$ to be a normal convex integrand is satisfied, in particular, if it is a piecewise continuous function of $t$, continuous in $(s,i)$ and convex in $u$. 
}
\end{remark}

\section{The case of a linear cost}\label{Section3}

\begin{problem}\label{IHProb}
Given  the initial data $(s_0,i_0)\in\mathcal{B}$  in the feasible region, minimize,  over all admissible controls $u\in U$,
\begin{itemize}
\item the {\em cost functional}
\begin{equation}\label{cf31}
J^\infty(u,s,i)=\int_0^{\infty}\big[\lambda_1 u(t)+ \lambda_2i(t)\big]\,dt,
\end{equation}
with $\lambda_2\ge0$, and $\lambda_1>0$,
\item under the {\em ICU constraint} on the trajectory of \eqref{stateeq}
\begin{equation}\label{icu}
i(t)\le i_M,\quad \forall t\in[0,+\infty).
\end{equation}
\end{itemize}
In the sequel, we refer to this formulation as problem $\mathcal{P}^\infty_{{\lambda_1},{\lambda_2}}$.
\end{problem}

\begin{remark}\label{ccwfc}{\em   
Non-emptiness of the set of viable controls with a corresponding finite cost is guaranteed by 
using, for example, the {\em greedy strategy}, $u_g$,
described in  Section \ref{higcs} 
till the herd immunity is reached, and $u_g=0$ afterwards.
Indeed,  as observed in Remark \ref{hirft}, with this strategy the herd immunity is reached in finite time. Hence, there exists $t_I\ge 0$ such that $s(t_I)\le\frac{\gamma}{\beta}$ and $u_g=0$ in $(t_I,+\infty)$. Moreover, since $s$ is strictly decreasing, there exists $\eps>0$ and $t_\eps\ge t_I$ such that
 $s(t)\le\frac{\gamma}{\beta}-\eps$ for every $t\in(t_\eps,+\infty)$ 
and hence
$$
\beta s(t)-\gamma\le-\eps\beta.
$$
Thus, for every $t\in(t_\eps,+\infty)$, we have
$$
i(t)=i_\eps\e^{\int_{t_\eps}^t[\beta s(\tau)-\gamma]\, d\tau}
\le i_I\e^{-\eps\beta(t-t_\eps)}.
$$
Then we have
$$
\int_{t_\eps}^{\infty}i(t)\,dt<+\infty,
$$
and this implies 
$$
J^\infty(u,s,i)=\int_0^{\infty}\big[\lambda_1 u(t)+ \lambda_2i(t)\big]\,dt
=\int_0^{t_\eps}\big[\lambda_1 u(t)+ \lambda_2i(t)\big]\,dt+\lambda_2\int_{t_\eps}^{\infty}i(t)\,dt<+\infty.
$$
}
\end{remark}

\begin{proposition}\label{PropV}
For every $\pr{s_0,i_0}\in\mathcal{B}$ with $i_0>0$ the problem $\mathcal{P}_{{\lambda_1},{\lambda_2}}^\infty$ has an optimal solution with a finite cost.
\end{proposition}

\proof Follows by Theorem \ref{gen_ex_th}, {\color{black}by taking into account that Remark \ref{ccwfc} ensures that there exists at least an admissible pair $(u,s,i)$ with a finite cost.}
\QED

Very important for our analysis is that, with an infinite time horizon, to reach herd immunity in finite time is a necessary condition for optimality. Indeed, the following consequence of Theorem \ref{uL1} holds true.

\begin{corollary}\label{corhnc}
If $u$ is an optimal control for problem $\mathcal{P}_{{\lambda_1},{\lambda_2}}^\infty$, then $s^u_\infty<\frac{\gamma}{\beta}$.
\end{corollary}

\proof
By Theorem \ref{uL1}, it is enough to observe that $u\in L^1(I)$. In fact, this is a consequence of the fact that optimal controls have a finite cost  (see Proposition \ref{PropV}), and of the inequality
$$
\int_0^\infty  u\,dt\le \frac{1}{\lambda_1}J^\infty(u,s,i)<\infty.
$$

\QED

\vspace{2ex}

\section{Reduction to finite horizon problems}\label{RFHP}

Besides the elements already discussed in the previous section (Corollary \ref{corhnc}), other optimality conditions can be stated by using Pontryagin's theorem on finite horizon problems suitably related to the (infinite horizon) original problem.

{\color{black}
\begin{problem}\label{ProbFixedtf_0_v2}
Let $T>0$. Given  the initial data $(s_0,i_0)\in\mathcal{B}$  in the feasible region, minimize,  over all admissible controls $u\in {L}^\infty\big((0,T);[0,u_{\max}]\big)$,
\begin{itemize}
\item the {\em cost functional}
\begin{equation}\label{cf31_v2}
J^T(u,s,i)=\int_0^{T}\big[ \lambda_2i(t)+\lambda_1 u(t)\big]\,dt,
\end{equation}
with $\lambda_2\ge0$ and $\lambda_1>0$,
\item under the {\em ICU constraint} on the trajectory of \eqref{stateeq}
\begin{equation}\label{icu_v2}
i(t)\le i_M\quad \forall t\in[0,T].
\end{equation}
\end{itemize}
In the sequel, we refer to this formulation as problem $\mathcal{P}^T_{{\lambda_1},{\lambda_2}}$.
\end{problem}
Let us remark that when $T=\infty$ we obtain the infinite horizon problem $\mathcal{P}^\infty_{{\lambda_1},{\lambda_2}}$ already introduced as Problem \ref{IHProb}.
}

Assume $u\in U$ to be an optimal control for problem $\mathcal{P}^\infty_{{\lambda_1},{\lambda_2}}$. By Corollary \ref{corhnc}, we have
$s_\infty^u<\frac{\gamma}{\beta}$, hence  the time $T^u_{\frac{\gamma}{\beta}}$, to reach the herd immunity level $s^u=\frac{\gamma}{\beta}$ with the control $u$, turns out to be finite. If $\lambda_2=0$, after time $T^u_{\frac{\gamma}{\beta}}$ the optimal control is identically $0$, because $i$ is decreasing and the ICU constraint is trivially satisfied, and any other choice of $u$ would lead to a bigger cost. This allows to use the finite horizon problem on $[0,T]$ to state other necessary (Pontryagin) optimality conditions in the case $\lambda_2=0$. Unfortunately, the same argument cannot be used if $\lambda_2>0$. Nevertheless, it will be shown that a reduction to finite horizon problems can be 
obtained by a $\Gamma$-convergence argument, covering also the case $\lambda_2>0$. 

\subsection{The case of a cost depending only on the control ($\lambda_2=0$)}

In the case $\lambda_2=0$ the following reduction theorem holds.

\begin{theorem}[Reduction to a finite horizon, case $\boldsymbol{\lambda_2=0}$]\label{i2fh}
Let $u\in U$ be an optimal control for problem $\mathcal{P}^\infty_{{\lambda_1},{\lambda_2}}$. 
If $\lambda_2=0$ then, for every $T> T^u_{\frac{\gamma}{\beta}}$, the restriction $u_|:=u_{|[0,T]}$ is an optimal control for the finite horizon problem $\mathcal{P}^T_{{\lambda_1},{\lambda_2}}$ with target objective $s(T)<\frac{\gamma}{\beta}$. 
\end{theorem}

\proof  By contradiction, assume that $u_|$ is not optimal. Then there exists $\ub\in L^\infty((0,T);[0,u_{\max}])$ with $s^\ub(T)<\frac{\gamma}{\beta}$ and such that
$$
J^T(\ub,s^\ub,i^\ub)<J^T(u_|,s^{u_|},i^{u_|}).
$$ 
Denoting by $\ub_0$ the extension of $\ub$ to the interval $(0,+\infty)$ by setting it to zero after $T$, using the local property of the systems of differential equations, we would have 
$i^\ub=i^{\ub_0}_{|[0,T]}$, and $s^\ub=s^{\ub_0}_{|[0,T]}$. Then, in particular, $s^{\ub_0}(T)=s^\ub(T)<\frac{\gamma}{\beta}$.
 This implies
\begin{equation*}
J^\infty(\ub_0,s^{\ub_0},i^{\ub_0})=
\int_0^T\lambda_1 \ub\,dt=J^T(\ub,s^\ub,i^\ub)<J^T(u_|,s^{u_|},i^{u_|})\le J^\infty(u,s^u,i^u),
\end{equation*}
which contradicts the optimality of $u$. \QED

\

Since, by Corollary \ref{corhnc}, $T^u_{\frac{\gamma}{\beta}}$ is finite, then 
 the infinite horizon  optimal control problem is brought back to a finite horizon OCP with a target objective. The latter has been studied in \cite{AFG2022,AFG2022C} and the structure of the (unique) optimal control has been characterized. According to Theorem 5.6 of \cite{AFG2022}, the greedy strategy is optimal also in the case of an infinite horizon (and without any target objective). 
This explains also what has been observed by Molina \& Rapaport, in the introduction to \cite{MR2022}, when they compare their result with those obtained in \cite{AFG2022}.

\subsection{The general case}

When the cost depends also from $i$, which, differently from the control, vanishes only as $t\to\infty$, the argument used before does not work. Nevertheless, in this section we show that (even when $\lambda_2>0$) a reduction to finite horizon problems can be  obtained by a $\Gamma$-convergence argument. 
$\Gamma$-convergence is a general notion of variational convergence introduced in 1976 by De Giorgi and Franzoni \cite{DGF75,DGF79}; we refer the unaccustomed reader also to the books of Dal Maso \cite{DMas} and Braides \cite{Bra02book}  for general properties and selected applications. {\color{black}To shorten notation, throughout this section all problems refer to the same fixed intial data $(s_0,i_0)$ in the viable set ${\mathcal B}$.}

%
%

\begin{theorem}\label{GCTconv}
For every increasing sequence of positive numbers $T_n\to\infty$, the sequence of problems $\mathcal{P}^{T_n}_{{\lambda_1},{\lambda_2}}$ $\Gamma$-converges to the limit problem 
$\mathcal{P}^\infty_{{\lambda_1},{\lambda_2}}$ in the following sense:
\begin{enumerate}
\item (liminf inequality) for every sequence $u_n\in U$ with $i^{u_n}\le i_M$, $u_n=0$ on $(T_n,+\infty)$, 
 and 
$u_n\weaks u$ in $L^\infty(I)$ we have $i^u\le i_M$ and 
$$
\liminf_{n\to\infty}J^{T_n}(u_n,s^{u_n},i^{u_n})\ge J^\infty(u,s^u,i^u);
$$ 
\item (recovery sequence) for every $u\in U$ with $J^\infty(u,s^u,i^u)<\infty$ {\color{black}with $i^u\le i_M$,} 
there exists $u_n\in U$ with $i^{u_n}\le i_M$, $u_n=0$ on $(T_n,+\infty)$, 
 and 
$u_n\weaks u$ in $L^\infty(I)$ such that
$$
\lim_{n\to\infty}J^{T_n}(u_n,s^{u_n},i^{u_n})= J^\infty(u,s^u,i^u).
$$ 
\end{enumerate}
\end{theorem}

\begin{remark}\label{rem_gamma}{\em 
\begin{enumerate}
\item The statement of Theorem \ref{GCTconv} can also be written in terms of usual $\Gamma$-convergence 
as $T\to\infty$ with respect to the weak* convergence on the space $L^\infty(I;[0,u_{\max}])\times W^{1,\infty}(I;\R^2)$
of the family of functionals
$$
F^T(u,s,i):=J^T(u,s,i)+\chi_{\Lambda}(u,s,i)+\chi_{i\le i_M}(i)+\chi_{u=0\mbox{ in }(T,+\infty)}(u)
$$
to the $\Gamma$-limit functional 
$$
F^\infty(u,s,i):=J^\infty(u,s,i)+\chi_{\Lambda}(u,s,i)+\chi_{i\le i_M}(i).
$$
\item The condition $u_n=0$ in $(T_n,+\infty)$ does not affect the problem  $\mathcal{P}^{T_n}_{{\lambda_1},{\lambda_2}}$; indeed, after $T_n$ the control can be set at any desired value without changing the value of the optimal control problem on the finite horizon $[0,T_n]$. 
\end{enumerate}
}
\end{remark}

{\color{black}

\begin{lemma}\label{lem_conv}
If $u_n\weaks u$ in $L^\infty(I;[0,u_{\max}])$, then $s^{u_n}\to s^u$ and $i^{u_n}\to i^u$ uniformly on every bounded subinterval $J\subset I$.
\end{lemma}

}

\proof   Since $u_n$, $s^{u_n}$ and $i^{u_n}$ are uniformly bounded, the derivatives of $s^{u_n}$ and  $i^{u_n}$ are uniformly bounded as well. Hence, $\|s^{u_n}\|_{1,\infty}$ and $\|i^{u_n}\|_{1,\infty}$ are bounded. {\color{black}Then the sequence $(s^{u_n},i^{u_n})$ is contained in a ball $B$ of the space $W^{1,\infty}(I;\R^2)$} on which (by the separability of $L^1$) the weak* topology is metrizable. 
Since bounded sets are weakly* relatively compact, then any subsequence of {\color{black} $(s^{u_n},i^{u_n})$} admits a weakly* converging subsequence. {\color{black} The application of Rellich compactness theorem ensures that weakly* converging sequences in $W^{1,\infty}(I;\R^2)$ are, up to subsequences, uniformly converging on every bounded subset $J\subset I$ (see for instance \cite[Theorem 8.8 and Remark 10]{Brezis2011}).
 We are then allowed to pass} to the limit in the state equations and, using the uniqueness of the solution, we get that these subsequence converges to {\color{black}$(s^{u},i^{u})$}. Since the limit is always the same and we are inside a metrizable set,
then we can conclude that the whole sequence converges to these limit.
Hence,  
$(s^{u_n},i^{u_n})\weaks (s^u,i^u)$ in $W^{1,\infty}(I;\R^2)$. The claimed uniform convergence on bounded subintervals comes again form the application of Rellich's theorem. \QED

\

\proof (of Theorem \ref{GCTconv}). {\em 1.}  By the previous lemma, the convergence assumption $u_n\weaks u$ in $L^\infty(I;[0,u_{\max}])$ implies 
that $i^{u_n}$ converges to $i^{u}$ uniformly {\color{black}on the bounded subintervals of $I$}, and this implies $i^u\le i_M$. 
%
{\color{black}Let us now remark that 
$$
\liminf_{n\to\infty}\int_0^{T_n}i^{u_n}=\liminf_{n\to\infty}\int_0^{\infty}i^{u_n}1_{[0,T_n]}\,dt\ge \int_0^{\infty}i^{u}\,dt
$$
simply by using the fact that the sequence $i^{u_n}1_{[0,T_n]}$ pointwisely converges to $i^{u}$ on $(0,+\infty)$ and Fatou's lemma (see, for instance, \cite[Theorem 11.31]{RudinPMA}).

By using the fact that $u_n=0$ on $(T_n,+\infty)$, 
the semicontinuity part {\em a.}\  of Theorem \ref{gen_ex_th} applied to the function $f_0(t,u,s,i)= u$, and the last inequality, we get
\begin{equation}\label{liJ}\begin{split}
\liminf_{n\to\infty}J^{T_n}(u_n,s^{u_n},i^{u_n})=&\liminf_{n\to\infty}\int_0^{T_n}\lambda_1 u_n+\lambda_2 i^{u_n}\,dt\\
=&\lambda_1 \liminf_{n\to\infty}\int_0^\infty u_n\,dt+\lambda_2\liminf_{n\to\infty}\int_0^{T_n}  i^{u_n}\,dt \\
\ge&\int_0^\infty\lambda_1 u+\lambda_2 i^{u}\,dt= J^\infty(u,s^u,i^u).
\end{split}\end{equation}}

\

{\em2.}  It remains to prove the existence of a recovery sequence. It can be constructed as follows.
Since $J^\infty(u,s^u,i^u)<\infty$,  the control $u$ belongs to $L^1(I)$ and, hence, $s_\infty^u<\frac{\gamma}{\beta}$, by Proposition \ref{uL1}.

Let us define
$$
u_n(t)=\begin{cases}
u(t)&\mbox{ if }t\in(0,T_n),\\
0&\mbox{ if }t\in(T_n,+\infty).
\end{cases}
$$
Since $T_n\to\infty$, for every $n$ large enough we have $s^{u_n}(T_n)=s^u(T_n)<\frac{\gamma}{\beta}$. This implies that $i^{u_n}\le i_M$. Indeed, before time $T_n$, 
we have $i^{u_n}=i^u\le i_M$, while for every time after the function $i^{u_n}(\cdot)$ is decreasing.

Moreover, $u_n\weaks u$ in $L^\infty(I)$ (by Lebesgue's theorem). 
Finally, by monotone convergence we have
\begin{eqnarray}\label{lsJ}
\lim_{n\to\infty}J^{T_n}(u_n,s^{u_n},i^{u_n})&=&\lim_{n\to\infty}\int_0^{T_n}\lambda_2 i^{u_n}+\lambda_1 u_n\,dt\\
&=&\int_0^\infty\lambda_2 i^{u}+\lambda_1 u\,dt= J^\infty(u,s^u,i^u)
\end{eqnarray}
and the theorem is completely proved. 
\QED

\

The next theorem is a consequence of the variational property of $\Gamma$-convergence.

\begin{theorem}\label{occonv} Let $u_T$ be an optimal control for $\mathcal{P}^{T}_{{\lambda_1},{\lambda_2}}$ with $T<\infty$, extended by $0$ in $(T,+\infty)$. Then 
\begin{enumerate}
\item there exists an increasing sequence $T_n\to\infty$ and $u\in U$ such that $u_{T_n}\weaks u$;
\item if $(T_n)$ is  sequence of positive numbers such that $T_n\to\infty$ and $u_{T_n}\weaks u$, then $u$ is an optimal control for problem $\mathcal{P}^{\infty}_{{\lambda_1},{\lambda_2}}$.
\end{enumerate}
\end{theorem}

\proof The first part of the statement immediately follows by the fact that $(u_T)$ is an equi-bounded family in $L^\infty(I)$ and the existence of a weak* converging sequence follows by Alaoglu's theorem. 

About assertion {\em2.}, by using Theorem \ref{esT}, it is easily seen
 that the optimal pair $(u_{T_n},s^{u_{T_n}},i^{u_{T_n}})$ weakly* converges in $U\times Y$  to $(u,s^{u},i^{u})$, and the claim follows by the variational property of $\Gamma$-convergence
(see, for instance, \cite[Corollary 7.17]{DMas}). \QED

\begin{remark}\label{gconvtec}{\em 
The $\Gamma$-convergence Theorem \ref{GCTconv} and its consequence Theorem \ref{occonv} allow us to deal with finite horizon problems. Indeed, suppose that $u_n$ is an optimal control for problem $\mathcal{P}^{T_n}_{{\lambda_1},{\lambda_2}}$ (which exists), extended by $0$ on $[T_n,+\infty)$.
Since $U$ is weakly* compact, then, up to a subsequence, $u_n\weaks u$ and $u$ is optimal for $\mathcal{P}^{\infty}_{{\lambda_1},{\lambda_2}}$.

Since $u$ is optimal, then $s^u_\infty<\frac{\gamma}{\beta}$ (see Corollary \ref{corhnc})  and, therefore, there exists $\eps>0$ and $T_\eps\in[0,+\infty)$ such that $s^u(t)<\frac{\gamma}{\beta}-\eps$ for every $t>T_\eps$.

On the other hand, {\color{black}by Lemma \ref{lem_conv},} $u_n\weaks u$ implies {\color{black}$s^{u_n}(T_\eps)\to s^u(T_\eps)$}. Hence, there exists $n_\eps\in\N$ such that 
$s^{u_n}(t)<\frac{\gamma}{\beta}$ for every $t>T_\eps$ and $n> n_\eps$. 
Since $T_n\to\infty$, then, possibly increasing $n_\eps$, we have $T_n>T_\eps$ for every $n> n_\eps$. Hence
 $s^{u_n}(T_n)<\frac{\gamma}{\beta}$ for every  $n> n_\eps$.

We are, essentially, saying that for any $T$ large enough, the optimal control $u_T$  of $\mathcal{P}^{T}_{{\lambda_1},{\lambda_2}}$ satisfies the final condition $s^{u_T}(T)<\frac{\gamma}{\beta}$.
We are, then, led to characterize the optimal controls for such $T$ large enough being allowed, in doing it, to use the final condition $s(T)<\frac{\gamma}{\beta}$ as a necessary condition (not as an assumption). {\color{black} As already explained in the introduction, this necessary optimality condition has a very important simplifying consequence in the application of Pontryagin's principle, because it implies normality (this will become clear in the proof of Theorem \ref{ThExtrA} and \ref{B0soc}).}}
%
\end{remark}
%
%

\section{The finite horizon problem: optimality conditions}\label{sec:fhphc}

\subsection{Pontryagin approach. General considerations}\label{Section4}

To write necessary conditions of optimality
let us introduce the adjoint variables $p_0\ge0$, $p_s,p_i\in\R$,
and the pre-Hamiltonian
$$
H(t,u,s,i,p_0,p_s,p_i)=p_0f_0(s,i,u)+p_s f_s+p_if_i,
$$
where $f_0(s,i,u)=\lambda_2i+\lambda_1u$ is the running cost function and  $f_s=-s(\beta-u)i$, $f_i=s(\beta-u)i-\gamma i$  are the dynamics of the state equations.
After some manipulations,  the pre-Hamiltonian turns out to be
$$
H(t,u,s,i,p_0,p_s,p_i)=p_0(\lambda_2i+\lambda_1u)+\eta s(\beta-u)i- \gamma p_i i,
$$
where $$\eta:=p_{i}-p_s.$$ 
In the sequel we use a constrained version of Pontryagin's theorem developed in \cite{BdlV2010,BdlVD2013}.
In particular, we refer to \cite{BdlVD2013} for the definition of the space of functions with bounded variation $BV([0,T])$ which is given by extending functions in a constant way on an open interval containing $[0,T]$.
 We adopt here also the notation
used in \cite{BdlV2010,BdlVD2013} of denoting the distributional derivative of a $BV$ function $f$ (which is a measure) by $df$,
instead than 
$f'$ that we reserve to measures which are absolutely continuous with respect to Lebesgue as, for instance, in the state equations.

\subsection{Optimality conditions for problem ${\mathcal{P}}^T_{\lambda_1,\lambda_2}$ with $T<\infty$}\label{ssOCP31}

By Pontryagin's theorem, given an optimal solution $(s,i,b)$, there exist a constant $p_0\in\{0,1\}$,  adjoint state real functions $p_s,\,p_i\in BV([0,T])$, a multiplier for the state constraint $\mu\in BV([0,T])$  with a nondecreasing representative (hence with measure  distributional derivative $d\mu\ge0$) such that $\mu(T^+)=0$ (recall that the functions are extended outside $[0,T]$), 
that satisfy the following properties.
\begin{itemize}
\item[(P1)] The non-degeneration property
\begin{equation}\label{ndc}
p_0+d\mu([0,T])>0.
\end{equation}
\item[(P2)] The complementarity condition
\begin{equation}\label{M2}
\int_{[0,{T}]}\big(i(t)-i_M)\ d\mu(t)=0.
\end{equation}
\item[(P3)] The conjugate equations with transversality conditions
\begin{equation*}
\begin{cases}
dp_s=-\eta (\beta-u) i,\\ dp_{i}=-\big(p_0\lambda_2+\eta (\beta-u)s-\gamma p_{i}\big)-d\mu,\\
p_s(T)=0,\ p_i(T^+)=0,
\end{cases}
\end{equation*}
which hold as equalities between measures on $[0,T]$.  We observe that the boundary condition for the costate $p_i$ is
given on the right limit in $T$, since $p_i$ could be discontinuous in $T$  if the measure $d\mu$ charges this point. On the contrary,
$p_s$ is continuous in $[0,T]$ since the derivative is absolutely continuous with respect to the Lebesgue measure.
\item[(P4)] The minimality property
$$
H\big(t,u(t),s(t),i(t),p_0,p_s(t),p_i(t)\big)=\inf_{u\in [0,u_{\max}]}\!H\big(t,u,s(t),i(t),p_0,p_s(t),p_i(t)\big),
$$
for almost every $t\in[0,T]$.
\item[(P5)] The conservation property
$$
H\big(t,u(t),s(t),i(t),p_0,p_s(t),p_i(t)\big)=k,
$$
with $k$ constant, for a.e.\  $t\in[0,T]$ (see \cite[Lemma 7.7]{ocbook}).
\end{itemize}

\subsection{Remarks and consequences}\label{cpc}

Throughout all this section we assume that $(s=s^{s_0,i_0,u},i=i^{s_0,i_0,u},u)$ be a solution of the control Problem \ref{ProbFixedtf_0_v2} with initial data $(s_0,i_0)\in\mathcal{B}$ and write consequences of Pontryagin's necessary conditions. Here we assume, moreover, that 
$s(T)<\frac{\gamma}{\beta}$; for our problem, this is not a restriction as observed
in Remark \ref{gconvtec}.

We have the following consequences.
\begin{itemize}
\item[(C1)] By (P3) {\color{black}and the definition of $\eta$} the jump condition $[\eta(t)]=[p_i(t)]=-[\mu(t)]\le 0$ holds for every $t\in[0,T]$ ({\color{black}where the jump of a BV function $f$ in a point $t$ is defined by $[f(t)]:=\dis\lim_{s\to t^+}f(s)-\dis\lim_{s\to t^-}f(s)=df(\{t\})$} and the inequality follows by the fact that $\mu$ is non-decreasing).
\item[(C2)]
By the conservation of the Hamiltonian
\begin{equation}\label{H=kf0}
p_0\big(\lambda_2i(t)+\lambda_1u(t)\big)+\eta(t) \big(\beta-u(t)\big)s(t)i(t)- \gamma p_i(t) i(t)=k,\mbox{ for a.e.\  }t\in[0,T].
\end{equation}
\item[(C3)]  (P3) implies that
\begin{equation}\label{deta}
d\eta=dp_i-dp_s=\eta (\beta-u)(i-s)-p_0\lambda_2+\gamma p_{i}-d\mu,
\end{equation}
or, owing to \eqref{H=kf0},
\begin{equation}\label{dechc}
d\eta=\eta (\beta-u)i+\frac{p_0\lambda_1u-k}{i}-d\mu.
\end{equation}
\item[(C4)] Since $s(T)<\frac{\gamma}{\beta}$ and $s$ is continuous, we have that $s<\frac{\gamma}{\beta}$ in an interval $(t_I,T]$ with $t_I<T$.
 This implies that $i'=\big(s(\beta-u) -\gamma\big) i<-\frac{\gamma}{\beta}ui\le0$, hence $i$ is decreasing and
	therefore $i<i_M$ in $(t_I,T]$. Then, by complementarity, we have $d\mu((t_I,T])=0$ which implies that $p_i$, and thus $\eta$ and $\psi$, are continuous in $(t_I,T]$. This implies $p_i(T)=0$ and $\eta(T)=0$.
\end{itemize}
Since the cost is linear in the control $u$, the minimum value of the Hamiltonian on
$K=[0,u_{\max}]$ is achieved when $u\in\{0,u_{\max}\}$.
Hence, setting the {\em switching function}
$$
\psi:=\eta si,
$$
the optimal control has to satisfy
\begin{equation}\label{cuSIRlinu}
b(t)=\begin{cases}
0,&\textnormal{ if }\psi(t)< p_0\lambda_1,\\[0ex]
u_{\max},&\textnormal{ if }\psi(t)>p_0\lambda_1,
\end{cases}
\end{equation}
for almost every $t\in[0,T]$, where $\psi$ denotes any pointwise representative of the switching function.
{\color{black} By $(C_4)$, we note that $\psi(T)=0$.}

The criterium \eqref{cuSIRlinu} allows to compute the value of $k$ in  \eqref{H=kf0}.

\begin{proposition}\label{etage0a0af} 
In (C2) we have $k=p_0\lambda_2i(T)$.
\end{proposition}

\proof 	
By having (C4) in mind, let us distinguish two cases.
\begin{enumerate}
\item If  $p_0=1$ then $\psi(T)=0<p_0\lambda_1$ and, by continuity, we have $\psi(t)<p_0\lambda_1$ and hence $u(t)=0$ a.e.\ in a left neighborhood of $T$ which we still call $(t_I,T]$;  \eqref{H=kf0} implies
$$
k=p_0\lambda_2 i(t)+\eta(t)\beta s(t) i(t)-\gamma p_i(t)i(t)\ \mbox{a.e.\ }t\in (t_I,T],
$$
and, by taking the limit as $t\rightarrow T-$, and invoking (C4), we obtain
$$
k=p_0\lambda_2 i(T).
$$
\item  If $p_0=0$ then \eqref{H=kf0} implies
\begin{equation}\label{kp00}
k=\eta(t)\big(\beta-u(t)\big) s(t) i(t)-\gamma p_i(t)i(t)\ \mbox{a.e.\ }t\in (t_I,T].
\end{equation}
By taking a sequence $t_n\to T$ on which the equality \eqref{kp00} holds and passing to the limit as $n\to\infty$
then we get $k=0$, which proves the proposition also in this case.
\end{enumerate}
\QED

\begin{proposition}\label{etaposg}  $\eta(t)\ge0$ for almost every $t\in[0,T]$.
\end{proposition}

\proof Let us recall that, by (C4), $\eta$ is continuous in a left neighborhood of $T$ and $\eta(T)=0$.
 
We note that, since $\mu$ is nondecreasing, $d\mu\ge0$ and, therefore, \eqref{dechc} implies
\begin{equation}\label{detai>0}
d\eta\le\eta (\beta-u)i+\frac{p_0\lambda_1u -k}{i}.
\end{equation}
Without loss of generality,  we assume $\eta$, and thus $\psi$, to be right-continuous.
Moreover
 $\eta$, and thus $\psi$, cannot have increasing jumps (by $(C1)$).

Arguing by contradiction, let us assume that there exists $t\in(0,t_0)$ such that $\eta(t)=\eta(t^+)<0$.
Since the switching function $\psi=\eta s i$ has the same sign as $\eta$, 
and since $p_0\ge0$,  
we have $\psi(t^+)<p_0\lambda_1$. On the other hand, since $\psi$ is right-continuous, one is able to find some
$\eps>0$ such that  $\psi<p_0\lambda_1$,
and hence $u=0$, a.e.\  in $J:=(t,t+\eps)$.   Owing to \eqref{detai>0} with $u=0$,
in $J$ we have
\begin{equation*}
d\eta\le\eta \beta i-\frac{k}{i}.
\end{equation*}
Since $k\ge 0$ (see Proposition \ref{etage0a0af}) , we have that
$$
d\eta\le \eta \beta i
$$
and $\eta \beta i<0$ in $J$, and therefore $\eta$ is negative and strictly decreasing in $J$.

As we shall see, this implies that $\eta<0$, and strictly decreasing, in $(t,T)$, which implies $\eta(T)<0$, 
thus contradicting the fact that  $\eta(T)=0$.

To prove the claim that $\eta<0$ in $(t,T)$, let us set
$$
t_1:=\sup\{s\in(t,T]\ :\ \eta(s)<0\}.
$$
(recall that $\eta$ is assumed to be right-continuous).
Assume, again by contradiction, that $t_1<T$. Since $\eta\in BV$, then there exists the left limit $\eta(t_1^-)\le0$. Since $\eta$ has no positive jumps in $[0,T]$,  it follows that $\eta(t_1^+)\le0$.  On the other hand,  $\eta(t_1^+)<0$ would contradict the definition of $t_1$. It follows that $\eta(t_1)=0$.
Nevertheless, in the interval $(t,t_1)$ we have that $\eta<0$, hence $\psi<0$.   Then, in $(t,T)$ we have that  $\psi<p_0\lambda_1$, hence $u=0$ a.e., hence $d\eta\le \eta (\beta-u)i<0$, hence $\eta$ is decreasing, which contradicts $\eta(t_1)=0$.  Our assumption on $t_1$ is wrong such that $t_1=T$ and $\eta<0$ in $(t,T)$. As we have already indicated, this contradicts $\eta(T)=0$ and, hence, our Proposition follows. \QED

\begin{corollary}\label{etage0a0} 
\begin{enumerate}
\item $\psi$ is non-negative almost everywhere in $(0,T)$;
\item $p_s$ is continuous, non-increasing and non-negative in $[0,T]$.
\end{enumerate}
\end{corollary}

\proof It is a straightforward consequence of the non-negativity of $\eta$, of the definition of $\psi$, of the first adjoint equation and the final condition $p_s(T)=0$.
\QED

Using \eqref{deta},  one easily proves the following result.

\begin{proposition}\label{prdpsi} The distributional derivative of $\psi$ is the measure given by
$$
d\psi=si(\gamma p_s-p_0\lambda_2)-si\,d\mu.
$$
\end{proposition}

The next theorem provides a characterization of singular arcs, whenever they occur.

\begin{theorem}\label{sing} 
Let  $(t_1,t_2)\subset\pp{0,T}$ be an interval in which  $\psi=p_0\lambda_1$. Then
the following hold true
\begin{enumerate}
\item if $p_0=1$, then
$\gamma p_s(t)-\lambda_2>0$ and $i(t)=i_M$, for every $t\in(t_1,t_2)$;
\item if $p_0=0$, then 
$\eta=p_i=p_s=0,$
in $(t_1,T]$ and $d\mu([0,t_1])>0$.
\end{enumerate}
\end{theorem}

\proof Let us consider the case when $p_0=1$ and prove the first assertion.
First of all,  since $\psi=\eta s i=p_0\lambda_1>0$, then we have $\eta>0$ in $(t_1,t_2)$.  Since $\psi$ is constant, on this interval, it follows that $d\psi=0$ in $(t_1,t_2)$.  Since $s$ and $i$ are strictly positive,  by Proposition \eqref{prdpsi},  we have
\begin{equation}\label{dmugps}
d\mu=\gamma p_s-\lambda_2, \mbox{ in }(t_1,t_2).
\end{equation}
Then, using the complementarity condition (P2)  and $\gamma>0$, we have that
\begin{equation}\label{psxM}
 \int_{(t_1,t_2)} \big(\gamma p_s(t)-\lambda_2\big)\big(i_M-i(t)\big)\,dt=0.
\end{equation}
Since $i$ is continuous, in order to prove that $i=i_M$,  it suffices to prove
that  $\gamma p_s-\lambda_2>0$.

Since $d\mu\ge0$, from  \eqref{dmugps} we have
\begin{equation}\label{psgen}
\gamma p_s-\lambda_2\ge0 \mbox{ in }(t_1,t_2).
\end{equation}
Now we prove that the strict inequality holds true.\\
Suppose now, by contradiction, that there exists $t_0\in(t_1,t_2)$ such that $\gamma p_s(t_0)-\lambda_2=0$.
By the adjoint equation $dp_s=-\eta (\beta-u)i$ we have that $dp_s\le0$ in $(t_1,t_2)$, hence $\gamma p_s-\lambda_2\le0$ for every
$t\in(t_0,t_2)$ and, by \eqref{psgen}, the equality holds in this interval. Then in $(t_0,t_2)$ we would have $dp_s=0$, hence $\eta=0$ by the adjoint equation, and, finally $\psi=0<p_0\lambda_1$, thus providing a contradiction. The first assertion is now completely proved.

\vspace{1ex}

Let us assume $p_0=0$ and prove {\sl 2.} 
In the interval  $(t_1,t_2)$, we have that $\psi=0$,  hence $\eta=0$;
 by the adjoint equations, we also have $dp_s=0$, hence there exists a constant $c$ such that $p_s=c$ on this interval.\\
Since $0=\eta=p_i-p_s$,  we also get $p_i=c$ on the aforementioned interval. By conservation of the Hamiltonian, and since $p_0=0$ implies $k=0$ (see Proposition \ref{etage0a0af}),
we have
\begin{equation}\label{chgci}
-\gamma c i(t)=0,
\end{equation}
hence $c=0$.   
By Corollary \eqref{etage0a0},  
 $p_s$ is nonincreasing. Since $p_s$ is zero on $(t_1,t_2)$ and at the end point $T$,  we must have $p_s=0$ in $(t_1,T]$.  By the first adjoint equation,  it follows that $\eta=0$ on $(t_1,T]$ which implies that $p_i=0$ on $(t_1,T]$. By \eqref{dechc}, we have $d\mu((t_1,T])=-d\eta((t_1,T])=0$.  Finally, the non degeneration condition \eqref{ndc} requires $d\mu\pr{\pp{0,T}}>0$, which implies $d\mu([0,t_1])>0$.

\QED

\begin{proposition}\label{x=x_M} 
If $i=i_c\in\left(0,i_M\right]$ is constant on an interval $(t_1,t_2)$, then there exists a positive constant $k_s$ such that
\begin{equation}\label{ocx=xM}
u(t)=\beta-\frac{\gamma}{s(t)}=\beta-\frac{\gamma}{k_s-\gamma i_ct},
\end{equation}
a.e.\ in the interval. The constant is given by
$$
k_s=s(t_1)+\gamma i_ct_1=s(t_2)+\gamma i_ct_2.
$$
Moreover, since $u\in\pp{0,u_{\max}}$, we have
$$
s(t_1)\le\frac{\gamma}{\beta-u_{\max}}\ \mbox{ and }\ s(t_2)\ge\frac{\gamma}{\beta}.
$$

\end{proposition}

\proof By the second state equation with $i=i_c$ we immediately have $s(\beta-u)=\gamma$.
Then the first becomes $s'=-\gamma i_c$. Integrating, we obtain  that there exists a constant $k_s$ such that
$s(t)=k_s-\gamma i_0 t$.  The moreover part of the statement follows by imposing $0\le\beta-\frac{\gamma}{s}\le u_{\max}$
and using the fact that $s$ decreases (recall that $i_c$ cannot be $0$).

\QED

\section{Structure of the optimal controls}\label{SectionExtremalsI}

Given an initial epidemic state  $(s_0,i_0)\in{\mathcal B}$ and a control $u\in U$, 
the ICU {\em saturation time} is defined as 
\[\tau_1^{s_0,i_0,u}:=\sup\big\{t\in (0,T]\ :\ i^{s_0,i_0,u}<i_M \mbox{ in }[0,t)\big\},\]
whenever the set on the right hand side is nonempty and $\tau_1^{s_0,i_0,u}:=0$ if it is empty. As usual, it will simply be denoted by $\tau_1$ when the initial conditions and the control can be easily deduced by the context.

Since the case $\lambda_2=0$ has already been completely studied in \cite{AFG2022}, we consider here only the case $\lambda_2>0$.
Then, for the whole section, and unless otherwise stated, all statements are meant to hold under the assumptions $\lambda_2>0$ and $\lambda_1>0$.  Moreover, to avoid trivialities, we always assume $s_0>0$ and $i_0>0$.

Let us start by characterizing the structure of the optimal controls for the problem ${\mathcal{P}}^T_{\lambda_1,\lambda_2}$, with $T<\infty$, in the case in which the ICU capacity is never saturated, that is, when $\tau_1=T$. Of course, this case always occurs if $i_M=1$, that is if an ICU constraint is not prescribed. 

\begin{theorem}\label{ThExtrA} Let us assume that 
$T<\infty$.
 If $\pr{s_0,i_0}\in\mathcal{B}$, $u$ is an optimal control of ${\mathcal{P}}^T_{\lambda_1,\lambda_2}$ with $s(T)<\frac{\gamma}{\beta}$ and 
$\tau_1^{s_0,i_0,u}=T$, then  
$u$ has a bang-bang structure $0-u_{\max}-0$ with at most two switches; in particular, for a.e.\ $t\in[0,T]$ we have
\begin{equation}\label{oca}
u(t)=\begin{cases}
0&\mbox{ if }t\in[0,\sie_C\wedge\sie_S],\\
u_{\max}&\mbox{ if }t\in(\sie_C\wedge\sie_S,\sie_S),\\
0&\mbox{ if }t\in[\sie_S,T],\\
\end{cases}
\end{equation}
where
\begin{eqnarray}
&&\sie_C:=\sup\{t\in[0,T]\ :\ \psi<p_0\lambda_1\mbox{ in }[0,t)\}\mbox{ with }\sup\varnothing:=0,\label{sigmaC}\\
&&\sie_S:=\inf\{t\in[0,T]\ :\ \psi<p_0\lambda_1\mbox{ in }(t,T]\}\mbox{ with }\inf\varnothing:=T\label{tauI}.
\end{eqnarray}
\end{theorem}

\begin{remark}{\em 
Before starting the proof, let us remark that the subscripts $C$ and $S$ in $\sigma_C$ and $\sigma_S$ stay for {\em Confinement}
and {\em Safety}, respectively. Indeed, whenever $\sigma_C<\sigma_S$, confinement is needed starting at  time $\sigma_C$ and ending at the safety time $\sigma_S$. The case in which $\sigma_C\ge\sigma_S$, which happens if and only if  $\psi<p_0\lambda_1$ a.e.\ in $[0,T]$, is a degenerate case in which the optimal strategy is $u\equiv0$, that is, do nothing.}
\end{remark}

\proof 
Since $s_0$, $i_0$ and $T$ are fixed, to improve readability, the states $i^{s_0,i_0,u}$ and $s^{s_0,i_0,u}$, corresponding to the specific initial conditions $(s_0,i_0)$ and control $u$, will be denoted  by 
$i$ and $s$, respectively. Consequently, almost all references to $T$ are dropped in the proof, but restored in the statement.

The assumption $\tau_1=T$ implies $i<i_M$ on $[0,T)$. Let us note that the same holds on $[0,T]$. Indeed, 
by (C4), the assumptions $T<\infty$ and $s(T)<\frac{\gamma}{\beta}$ imply that $i$ is strictly decreasing in a left neighborhood of $T$.
This implies the following consequences. 
\begin{enumerate}
\item $d\mu([0,T])=0$, by complementarity.
\item $p_0=1$ by non-degeneracy.
\item $p_i$, $\eta$ and $\psi$ are continuous in $[0,T]$ (and $p_i(T)=0$, $\eta(T)=0$ and $\psi(T)=0$ as already observed in (C4)),
by {\em1.}\ and the adjoint equations.
\item  Singular arcs are excluded, that is $\psi=p_0\lambda_1$ in a subinterval $(t_1,t_2)\subset[0,T]$ is impossible. Indeed, 
by Theorem \ref{sing},  we would have $i=i_M$, which contradicts $i<i_M$ on $[0,T]$.
\item By the first item and Proposition \eqref{prdpsi}, it holds $d\psi=s i(\gamma p_s-p_0\lambda_2)$ and  we have 
$$
d\psi\ge0 \ \iff\ \gamma p_s-p_0\lambda_2\ge0\ \iff\ p_s\ge\frac{p_0\lambda_2}{\gamma}.
$$
\item $\psi(T)<p_0\lambda_1$. Follows by the facts that $\psi(T)=0$ and $p_0=1$. This implies that the set in the definition of  $\sigma_S$ is non-empty and $\sigma_S<T$.
\end{enumerate}
By $6.$ we have $\psi<p_0\lambda_1$ in $(\sigma_S,T)$ and, by \eqref{cuSIRlinu}, $u=0$ in this interval.

If $\sigma_S=0$ then we have $u=0$ in $(0,T)$ and \eqref{oca} is true because, by \eqref{tauI}, we have $\sigma_C\wedge\sigma_S=\sigma_S=0$.

Let us now assume that $\sigma_S>0$. Then, by continuity, $\psi(\sigma_S)=p_0\lambda_1$. 
By Corollary \ref{etage0a0} $p_s$ is non-increasing 
and, by 6., the derivative of $\psi$ can change sign only once in the interval $(0,\sigma_S)$.
Therefore, and since (by 5.) singular arcs cannot occur,  the superlevel set $\{\psi>p_0\lambda_1\}$ turns out to be the  interval, possibly empty, $(\sigma_C\wedge\sigma_S,\sigma_S)$ as defined in \eqref{tauI}, and the conclusion follows by \eqref{cuSIRlinu}. \QED

\begin{corollary}\label{corThExtrA} Let us assume that 
$T<\infty$ and $\pr{s_0,i_0}\in\mathcal{A}\setminus\partial\mathcal{A}$. If $u^T$ is an optimal control of ${\mathcal{P}}^T_{\lambda_1,\lambda_2}$ with $s(T)<\frac{\gamma}{\beta}$,
then $u^T$ has the bang-bang structure \eqref{oca}-\eqref{tauI}.
\end{corollary}

\proof  If $(s_0,i_0)\in\mathcal{A}\setminus\partial\mathcal{A}$, then we have 
\begin{equation*}\label{sICUe}
i<i_M\mbox{ on }[0,T].
\end{equation*}  
Indeed, in the sense of distributions, 
\[\dfrac{d}{dt}\big(i-\Phi_{0}(s)\big)=\big((\beta-u^T) s-\gamma\big)i+(\beta-u^T) s i\ \min\big\{\frac{\gamma}{\beta s}-1,0\big\}\le0.\]
As a consequence, $i(t)\leq \Phi_{0}\pr{s(t)}+i_0-\Phi_{0}(s_0)<\Phi_{0}\pr{s(t)}$, which implies the desired inequality.

Then, the claim follows by Theorem \ref{ThExtrA}.
\QED

\begin{remark}\label{remIM=1}{\em The case considered in Corollary \ref{corThExtrA} arises, for instance, if one set\dg{s} $i_M=1$, that is when the ICU constraint is not prescribed.}
\end{remark}

\begin{remark}{\em
The interval $(\sie_S,T)$ is always non-empty (by {6.}).
The interval $(\sie_C,\sie_S)$ may be empty or not, depending on the initial condition $(s_0,i_0)$ and the ratio between the coefficients $\lambda_2$ and $\lambda_1$. The numerical simulation performed in Example \ref{ita_cov_noICU} falls in this framework. 
}
\end{remark}

\noindent By taking to $\Gamma$-limit as $T\to\infty$ (see Theorems \ref{GCTconv} and \ref{occonv}) we get the following result about problem ${\mathcal{P}}^\infty_{\lambda_1,\lambda_2}$.

\begin{theorem}\label{IHA} Let us assume that 
$\pr{s_0,i_0}\in\mathcal{A}\setminus\partial\mathcal{A}$.
The problem ${\mathcal{P}}^\infty_{\lambda_1,\lambda_2}$ admits an optimal control $u^\infty$ with the following structure: there exist  $0\le\sigma_1\le\sigma_2<\infty$ such that 
\begin{equation}\label{ocaf}
u^\infty(t)=\begin{cases}
0&\mbox{ if }t\in(0,\sigma_1),\\
u_{\max}&\mbox{ if }t\in(\sigma_1,\sigma_2),\\
0&\mbox{ if }t\in(\sigma_2,\infty).\\
\end{cases}
\end{equation}
\end{theorem}

\proof  For $T<\infty$, let $u^T$ be an optimal control for problem ${\mathcal{P}}^T_{\lambda_1,\lambda_2}$. By Theorem \ref{ThExtrA},
there exist  $0\le\sie_1\le\sie_2<T$ such that 
\begin{equation*}
u^T(t)=\begin{cases}
0&\mbox{ if }t\in(0,\sie_1),\\
u_{\max}&\mbox{ if }t\in(\sie_1,\sie_2),\\
0&\mbox{ if }t\in(\sie_2,T).\\
\end{cases}
\end{equation*}
Then, there exist a sequence of positive numbers  $T_n\to\infty$ and two times $0\le\sigma_1^\infty\le\sigma_2^\infty$ such that 
$\sigma^{T_n}_1\to\sigma_1^\infty$ and $\sigma^{T_n}_2\to\sigma_2^\infty$.  It follows that 
$u^{T_n}\weaks u^\infty$ with 
\begin{equation*}
u^\infty(t)=\begin{cases}
0&\mbox{ if }t\in(0,\sigma_1^\infty),\\
u_{\max}&\mbox{ if }t\in(\sigma_1^\infty,\sigma_2^\infty),\\
0&\mbox{ if }t\in(\sigma_2^\infty,\infty).\\
\end{cases}
\end{equation*}
By Theorem \ref{occonv}, $u^\infty$ is an optimal control for ${\mathcal{P}}^\infty_{\lambda_1,\lambda_2}$.

If $\sigma_1^\infty<\sigma_2^\infty$  then we have $\sigma_2^\infty<\infty$. Indeed, otherwise the cost would be infinite. 
Then,  the theorem is proved in this case.

If $\sigma_1^\infty=\sigma_2^\infty$ then $u=0$ and the theorem is proved, with $\sigma_1=\sigma_2$, also in this case.

\noindent Then $u^\infty$ takes always the form \eqref{ocaf} and the theorem is completely proved.

\QED

\begin{theorem}\label{B0soc} Let $T<\infty$. Assume that $u$ be an optimal control for problem ${\mathcal{P}}^T_{\lambda_1,\lambda_2}$ with $s(T)<\frac{\gamma}{\beta}$. 
 If $(s_0,i_0)\in\mathcal{B}\setminus\mathcal{A}$ and $i_0<i_M$, then the control $u$ must have the following structure: 
\begin{equation}\label{uocbar1}
u(t)=\begin{cases}
0&\mbox{ if }t\in(0,\sigma_C\wedge\sigma_S),\\
u_{\max}&\mbox{ if }t\in(\sigma_C\wedge\sigma_S,\tau_1\wedge\sigma_S),\\
0&\mbox{ if }t\in(\tau_1\wedge\sigma_S,\tau_1),\\
\dis\beta-\frac{\gamma}{s(\tau_2)+\gamma i_M(\tau_2-t)}&\mbox{ if }t\in(\tau_1\wedge\sigma_S,\tau_2\wedge\sigma_S),\\[2ex]
u_{\max}&\mbox{ if }t\in(\tau_2\wedge\sigma_S,\sigma_{S}),\\
0&\mbox{ if }t\in(\sigma_S,T), 
\end{cases}
\end{equation}
where $\sigma_C$ and $\sigma_S$ are defined as in \eqref{sigmaC} and  \eqref{tauI}, while 
\[\tau_1:=\sup\big\{t\in (0,T]\ :\ i<i_M \mbox{ in }[0,t)\big\},\]
 and, if $\tau_1<T$, 
\begin{equation}\label{deftau2}
\tau_2:=\sup\big\{t\in[\tau_1,T]\ :\ i=i_M \textnormal{ in }[\tau_1,t]\big\}.
\end{equation}
\end{theorem}

\begin{remark}{\em Note that some among the time intervals in \eqref{uocbar1} might be empty.
For instance, $(\tau_1\wedge\sigma_S,\tau_1)=\varnothing$ if $\sigma_S\ge\tau_1$.
}
\end{remark}

\proof If $\tau_1=T$ then the previous Theorem \ref{ThExtrA} applies and $u$ has the structure $0-u_{\max}-0$, as claimed.
Let us then  assume that $\tau_1<T$ (hence $i(\tau_1)=i_M$).

Let us start with some structural remarks on $[0,\tau_1)$.
The assumption $i_0<i_M$ implies $\tau_1>0$. Then the interval $[0,\tau_1)$ is nonempty  and, inside, we have $i<i_M$. 
The following consequences hold.
\begin{enumerate}
\item[(a)]$d\mu([0,\tau_1))=0$, by complementarity.
\item[(b)] $p_i$, $p_s$, $\eta$ and $\psi$, are continuous in  $[0,\tau_1)$, by the adjoint equations and the previuos item.
\item[(c)] Singular arcs cannot occur in $[0,\tau_1)$. Indeed,
\begin{itemize}
\item if $\psi=p_0\lambda_1=0$ in a subinterval $(t_1,t_2)\subset[0,\tau_1)$, then $p_0=0$ and, owing to Theorem \ref{sing}, we have  $d\mu([0,t_1])>0$ which is impossible since $[0,t_1]\subset [0,\tau_1)$ and
$d\mu([0,\tau_1))=0$ by (a);  
\item $\psi=p_0\lambda_1>0$ in a subinterval $(t_1,t_2)\subset[0,\tau_1)$ is impossible; indeed, this would imply $p_0=1$ and by Theorem \ref{sing}  we would have $i=i_M$, against the definition of $\tau_1$. 
\end{itemize}
\item[(d)] By the first assertion and Proposition \ref{prdpsi}, on $[0,\tau_1)$ we have
$d\psi=si(\gamma p_s-p_0\lambda_2)$;  hence 
$$
d\psi\ge0 \ \iff\ \gamma p_s-p_0\lambda_2\ge0\ \iff\ p_s\ge\frac{p_0\lambda_2}{\gamma},
$$
in $[0,\tau_1)$.
\item[(e)]  Starting from an initial condition $(s_0,i_0)\in\mathcal{B}\setminus\partial\mathcal{B}$  we can exclude that $u=u_{\max}$ in the whole interval $(0,\tau_1)$.

First of all, this can be excluded if $s_0\le\frac{\gamma}{\beta-u_{\max}}$, because  in such a case the control $u=u_{\max}$ would
keep the trajectory in the interior of $\mathcal{B}$ and could not reach the level $i_M$ at time $\tau_1$. Indeed, the function $i^{s_0,i_0,u_{\max}}(s)$ has a maximum in $s=\frac{\gamma}{\beta-u_{\max}}$ and $s$ is a decreasing function of time.

Also in the case $s_0>\frac{\gamma}{\beta-u_{\max}}$, the associated trajectory $(s^{s_0,i_0,u_{\max}},i^{s_0,i_0,u_{\max}})$ would be kept in the interior of $\mathcal{B}$, so contradicting $i^{s_0,i_0,u_{\max}}(\tau_1)=i_M$.  In fact, by the laminar flow property (see 1.-4.\ of Section \ref{EVA}) starting from a position $(s_0,i_0)$ under the curve $\Phi_{\max}$, the control $u_{\max}$ generates a trajectory that is always strictly under $\Phi_{\max}$ as long as $s\ge\frac{\gamma}{\beta-u_{\max}}$, where it takes its maximum value strictly smaller than $i_M$. Hence, it cannot reach the boundary at time $\tau_1$.
\item[(f)]  We have $p_0=1$. Indeed, assume by contradiction that $p_0=0$. By Corollary \ref{etage0a0}, 
 we have  $\psi\ge0=p_0\lambda_1$ and $p_s\ge0=\frac{p_0\lambda_2}{\gamma}$ in $[0,T]$ .
By (d) this implies that $\psi$ is non-decreasing in $[0,\tau_1)$. 
Since singular arcs are excluded, then we have $\psi>p_0\lambda_1$, hence $u=u_{\max}$ in $[0,\tau_1)$ which is excluded by (e).
\end{enumerate}

\noindent\emph{Step 1.} Let us prove that, in the interval  $\pr{0,\tau_1}$ the optimal control has a bang-bang structure $0-u_{\max}-0$ with at most two switches, that is
\begin{equation}\label{oca01}
u(t)=\begin{cases}
0&\mbox{ if }t\in(0,\sigma_C),\\
u_{\max}&\mbox{ if }t\in(\sigma_C,{\sigma}_{CE}),\\
0&\mbox{ if }t\in({\sigma}_{CE},\tau_1),
\end{cases}
\end{equation}
where  $\sigma_C$ has been defined defined in \eqref{sigmaC}, while
$$
\sigma_{CE}:=\sup\{t\in(\sigma_C,\tau_1)\ :\ \psi(t)>p_0\lambda_1\},\quad \sup\varnothing=\sigma_C,
$$
is the {\em Exiting} time.\\
Indeed, since $p_s$ is non-increasing (see Corollary \ref{etage0a0}), {\color{black} it can cross the level $\frac{p_0\lambda_2}{\gamma}$ at most one time. On the other hand, the sign of $d\psi$ is related to the value of $p_s$ by point $(d)$. It follows that}  $\psi$ is quasi-concave (in particular, it cannot increase after having started to decrease). Since, moreover, singular arcs cannot occur (by (c)),
then $\psi$ can cross the level $p_0\lambda_1$ at most twice in the interval $[0,\tau_1)$, that is in $\sigma_C$ if it is not zero and in $\sigma_{CE}$ if is not equal to $\tau_1$.  Thus, the claimed structure \eqref{oca01} follows by \eqref{cuSIRlinu}.

\

\noindent\emph{Step 2.} Since $i(\tau_1)=i_M$,  it follows that the time $\tau_2$ is well defined, that is, the set on the right hand side in \eqref{deftau2} is nonempty. By the assumption $s(T)<\frac{\gamma}{\beta}$, one has $\tau_2<T$ (see (C4), Section \ref{cpc}). We focus on the case when $\tau_1<\tau_2$ for which, in $(\tau_1,\tau_2)$, we have $i=i_M$. Then, the optimal control in $(\tau_1,\tau_2)$ is given by \eqref{ocx=xM}, that is
\begin{equation}\label{bpar}
u(t)=\beta-\frac{\gamma}{s(\tau_2)+\gamma i_M(\tau_2-t)}.
\end{equation}
{Since $u\in[0,u_{\max}]$,  we have (see Proposition \ref{x=x_M})
\begin{equation}\label{sbounds}
\frac{\gamma}{\beta}\le s(\tau_2)\le \frac{\gamma}{\beta-u_{\max}} .
\end{equation}}
Since $s$ is strictly decreasing,  it follows that
\begin{equation}\label{sgb}
s(t)>s(\tau_2)\ge \frac{\gamma}{\beta}\ \mbox{ for every }t\in[\tau_1,\tau_2).
\end{equation}
We distinguish the following two cases:
\begin{enumerate}
\item $s(\tau_2)=\frac{\gamma}{\beta}$. 
Since  $s$ is decreasing we have $s<\frac{\gamma}{\beta}$ in $(\tau_2,T]$. Since $i(\tau_2)=i_M$ and herd immunity has been reached, we have $i<i_M$ in $(\tau_2,T]$. Then we are under the assumptions of Theorem \ref{ThExtrA} with $s_0=\frac{\gamma}{\beta}$ 
and  the optimal control has the structure $0-u_{\max}-0$, with at most two switches, in the interval $(\tau_2,T]$.
\item $s(\tau_2)>\frac{\gamma}{\beta}$. We observe that 
$\psi(\tau_2^+)<p_0\lambda_1$ is excluded. Indeed, in such case there would be an interval $(\tau_2,\sigma)$ in which $u=0$ and the constraint $i=i_M$ would be violated.  Then we have $\psi(\tau_2^+)\ge p_0\lambda_1$. \\
Subcases:
\begin{enumerate}
\item $p_s(\tau_2)\le\frac{p_0\lambda_2}{\gamma}$. Since $p_s$ is non-increasing, by Proposition \ref{prdpsi} we have $d\psi\le0$,  on $(\tau_2,T]$, that is also $\psi$ is non-increasing. On the other hand, it cannot stay equal to $p_0\lambda_1$ since, otherwise, recalling that $p_0=1$, we would have $i=i_M$ in a right neighborhood of $\tau_2$, against its definition, then we have $\psi<p_0\lambda_1$, and, hence $u=0$, in $(\tau_2,T]$.
\item $p_s(\tau_2)>\frac{p_0\lambda_2}{\gamma}$. By continuity of $s$ and $p_s$, there exists $\sigma>\tau_2$ such that 
$p_s>\frac{p_0\lambda_2}{\gamma}$ and $s>\frac{\gamma}{\beta}$ in $(\tau_2,\sigma)$. Let us 
 we prove that
\begin{equation}\label{claimi<iM}
i<i_M\mbox{ in a right neighborhood of }\tau_2.
\end{equation}
We assume, by contradiction, the existence of a strictly decreasing sequence of points $t_n\in(\tau_2,\sigma)$ with $i(t_n)=i_M$ and such that $t_n\to\tau_2$.
We claim that this implies that
\begin{equation}\label{tnM}
i=i_M\mbox{ in }[t_n,t_{n+1}]\mbox{ for every $n$} .
\end{equation}
Again,  by contradiction, should there exist $\hat t_n\in(t_{n+1},t_{n})$ such that $i(\hat t_n)<i_M$, the continuity of $i$ implies the existence of
some open interval $\hat t_n\in(\underline{a},\overline{a})\subseteq[t_{n+1},t_{n}]$ such that $i<i_M$ inside and $i=i_M$ at the boundary (one may take, for instance, $\underline{a}=\inf\{t\in[t_{n+1},\hat t_{n}]\ :\ i(t)<i_M\}$ and $\overline{a}=\sup\{t\in[\hat t_n,t_{n}]\ :\ i(t)<i_M\}$).
By complementarity,  we have that $d\mu((\underline{a},\overline{a}))=0$ and, therefore, $d\psi=si(\gamma p_s-p_0\lambda_2)$ so that  $\psi$ is continuous in $(\underline{a},\overline{a})$. 
Let us consider the following cases.
\begin{enumerate}
\item $\psi(\underline{a}^+)< p_0\lambda_1$. Then $\psi< p_0\lambda_1$ in a right neighborhood of $\underline{a}$ (in which $u=0$),
but this is excluded because $i$ would be increasing, but $i(\underline{a})=i_M$ and the constraint $i\le i_M$ would be violated.
\item $\psi(\underline{a}^+)\ge p_0\lambda_1$. Since $p_s>\frac{p_0\lambda_2}{\gamma}$ in $(\underline{a},\overline{a})$ and $d\mu(\underline{a},\overline{a})=0$ then we would have $\psi> p_0\lambda_1$, hence $u=u_{\max}$, in $(\underline{a},\overline{a})$, which would imply that $i$ is strictly decreasing in the interval against the fact that $i=i_M$ at the boundary.
\end{enumerate}
Since we always obtain a contradiction, we have proved that \eqref{tnM} holds true. Since the sequence $t_n$ converges to $\tau_2$ then, by continuity, this implies that
$i=i_M$ on $(\tau_1,\tau_2+\eps)$ with $\eps>0$, against the definition of $\tau_2$. This proves \eqref{claimi<iM} and
there exists an interval $(\tau_2,\sigma)$ in which we have $i<i_M$.
\end{enumerate}
By complementarity,  we have $d\mu(\pr{\tau_2,\sigma})=0$. 
Since $\psi(\tau_2^+)\ge p_0\lambda_1$ and $p_s(\tau_2)>\frac{p_0\lambda_2}{\gamma}$, then we have 
$\psi> p_0\lambda_1$ in a right neighborhood of $\tau_2$ in which $u=u_{\max}$ and $i$ is decreasing (keep in mind the upper bound on $s(\tau_2)$ and, thus, on every $s(t)$, for $t\geq \tau_2$). 
Let us define 
$$
\sigma_{2E}:=\sup\{t\in(\tau_2,T]\ :\ \psi> p_0\lambda_1\mbox{ a.e.\ in }(\tau_2,t)\}.
$$
We have $\psi(\sigma_{2E}^+)\le p_0\lambda_1$ and $i(\sigma_{2E})<i_M$, and $p_s(\sigma_{2E})\le\frac{p_0\lambda_2}{\gamma}$
(because, otherwise, $\psi$ should increase against the definition of $\sigma_{2E}$). 
Since $p_s$ is non-increasing, by Proposition \ref{prdpsi}, after this time we have $d\psi\le0$, that is also $\psi$ is non-increasing. On the other hand, it cannot stay equal to $p_0\lambda_1$ because, otherwise, $i=i_M$ in a right neighborhood of $\sigma_{2E}$, which is impossible by continuity of $i$. Then we have $\psi<p_0\lambda_1$, and hence $u=0$, in $(\sigma_{2E},T]$. 
And, of course, we get $\sigma_{2E}=\sigma_S$.

Summarizing, in this second case $u$ must have the structure
\begin{equation*}
u=\begin{cases}
\umax&\mbox{ in }(\tau_2\wedge\sigma_S,\sigma_{S}),\\
0&\mbox{ in }(\sigma_S,T).
\end{cases}
\end{equation*}
\end{enumerate}
Putting together with the first case:
\begin{equation}\label{bt2T}
u=\begin{cases}
0&\mbox{ in }(\tau_2,\sigma_{2C}\wedge\sigma_S),\\ 
\umax&\mbox{ in }(\sigma_{2C}\wedge\sigma_S,\sigma_{S}),\\
0&\mbox{ in }(\sigma_S,T), 
\end{cases}
\end{equation}
where 
$$
\sigma_{2C}:=\sup\{t\in(\tau_2,T]\ :\ \psi(t)<p_0\lambda_1\mbox{ a.e.\ in }(\tau_2,t)\},\quad \sup\varnothing=\tau_2.
$$
{\em Step 3.} Summarizing, we have proved that $u$ has the following structure
\begin{equation}\label{uocbar1pre}
u(t)=\begin{cases}
0&\mbox{ if }t\in(0,\sigma_C\wedge\sigma_S),\\
\umax&\mbox{ if }t\in(\sigma_C\wedge\sigma_S,{\sigma}_{CE}\wedge\sigma_S),\\
0&\mbox{ if }t\in({\sigma}_{CE}\wedge\sigma_S,\tau_1\wedge\sigma_S),\\
\dis\beta-\frac{\gamma}{s(\tau_2)+\gamma i_M(\tau_2-t)}&\mbox{ if }t\in(\tau_1\wedge\sigma_S,\tau_2\wedge\sigma_S),\\[2ex]
0&\mbox{ in }(\tau_2\wedge\sigma_S,\sigma_{2C}\wedge\sigma_S),\\ 
\umax&\mbox{ in }(\sigma_{2C}\wedge\sigma_S,\sigma_{S}),\\
0&\mbox{ in }(\sigma_S,T), 
\end{cases}
\end{equation}
The final form \eqref{uocbar1} is achieved by taking into account the following remarks.   
\begin{enumerate}
\item If $\sigma_S>\tau_1$ and $\tau_2>\tau_1$, that is if the boundary arc is non-empty,
then the regime $\umax-0$ immediately before $\tau_1$ cannot occur.

This is because {\color{black} we are speaking of the interval $[0,\tau_1)$ considered in Step 1, characterized in \eqref{oca01} and for which the structural remarks (a)-(f) of the inital part of the proof hold.} {\color{black}Arguing by contradiction,} in order to pass from being larger than $p_0\lambda_1$ {\color{black}in the regime $u_{\max}$ (see \eqref{cuSIRlinu})}  to being smaller, $\psi$ has to decrease. Since $p_s$ is non-increasing in $[0,T]$ {\color{black}(Corollary \ref{etage0a0}, 2.)}, when $\psi$ has start to decrease, {\color{black}by (d)} we have $p_s<\frac{p_0\lambda_2}{\gamma}$. {\color{black}By monotonicity of $p_s$, this last inequality} holds, in particular, in the whole interval $({\sigma}_{CE}\wedge\sigma_S,T)$. Then, in this interval $\psi$ is decreasing {\color{black}(again by (d))} and could not reach the level $p_0\lambda_1$ {\color{black} necessary for the appearance} of the boundary arc {\color{black}(see \eqref{cuSIRlinu})} starting from being smaller (in the preceding stage $u=0$).

Then, if $\sigma_S\le\tau_1$ the regime $\umax-0$ can be replaced by $\umax$ and the optimal control can be written in the final form \eqref{uocbar1}{\color{black}, till to time $\tau_2\wedge\sigma_S$.}
\item Analogously, the regime $0-\umax$ cannot occur after $\tau_2$.
\item The stage $\umax$ after $\tau_2$ does not occur if $\lambda_2=0$ as proven in \cite{AFG2022}, but can occur if $\lambda_2>0$ (see figure \ref{ITA1}).
\end{enumerate}

\QED

\begin{theorem} Let us assume that 
$\pr{s_0,i_0}\in\mathcal{B}\setminus\mathcal{A}$ and $\tau_1<\infty$.
The problem ${\mathcal{P}}^\infty_{\lambda_1,\lambda_2}$ admits an optimal control $u^\infty$ with the following structure: there exist  $0\le\tau_0\le \tau_1\le\tau_2\le\tau_3<\infty$ such that 
\begin{equation}\label{ocIHBf}
u^\infty(t)=\begin{cases}
0&\mbox{ if }t\in(0,\tau_0),\\
\umax&\mbox{ if }t\in(\tau_0,\tau_1),\\
\dis\beta-\frac{\gamma}{s(\tau_2)+\gamma i_M(\tau_2-t)}&\mbox{ if }t\in(\tau_1,\tau_2),\\[2ex]
\umax&\mbox{ in }(\tau_2,\tau_3),\\
0&\mbox{ in }(\tau_3,\infty), 
\end{cases}
\end{equation}
\end{theorem}

\proof  For $T<\infty$, let $u^T$ be an optimal control for problem ${\mathcal{P}}^T_{\lambda_1,\lambda_2}$. By Theorem \ref{B0soc} ,
there exist  $0\le\sie_0\le\sie_1\le\tau_1^T\le \tau_2^T\le\sie_3<T$ such that 
\begin{equation*}
u^T(t)=\begin{cases}
0&\mbox{ if }t\in(0,\sie_0),\\
\umax&\mbox{ if }t\in(\sie_0,\sie_1),\\
0 &\mbox{ if }t\in(\sie_1,\tae_1),\\
\dis\beta-\frac{\gamma}{s(\tae_2)+\gamma i_M(\tae_2-t)}&\mbox{ if }t\in(\tae_1,\tae_2),\\[2ex]
\umax&\mbox{ in }(\tae_2,\sie_3),\\
0&\mbox{ in }(\sie_3,T), 
\end{cases}
\end{equation*}
Then, there exist a sequence of positive numbers  $T_n\to\infty$ and five times $0\le\sigma_0^\infty\le\sigma_1^\infty\le\tau_1^\infty\le \tau_2^\infty\le\sigma_3^\infty$ such that 
$\sigma^{T_n}_i\to\sigma_i^\infty$, $i=1,2,3$, and $\tau^{T_n}_j\to\tau_j^\infty$, $j=1,2$.  It follows that 
$u^{T_n}\weaks u^\infty$ with
\begin{equation}\label{uipre}
u^\infty(t)=\begin{cases}
0&\mbox{ if }t\in(0,\sigma_0^\infty),\\
\umax&\mbox{ if }t\in(\sigma_0^\infty,\sigma_1^\infty),\\
0 &\mbox{ if }t\in(\sigma_1^\infty,\tau_1^\infty),\\
\dis\beta-\frac{\gamma}{s(\tau_2^\infty)+\gamma i_M(\tau_2^\infty-t)}&\mbox{ if }t\in(\tau_1^\infty,\tau_2^\infty),\\[2ex]
\umax&\mbox{ in }(\tau_2^\infty,\sigma_3^\infty),\\
0&\mbox{ in }(\sigma_3^\infty,\infty), 
\end{cases}
\end{equation}
By Theorem \ref{occonv}, $u^\infty$ is an optimal control for ${\mathcal{P}}^\infty_{\lambda_1,\lambda_2}$.

If $\tau_1^\infty=\infty$ then $u^\infty=0-u_{\max}-0$ and the claim is true with $\tau_1=\tau_2$.

If $\tau_1^\infty<\infty$ we have the following subcases.
\begin{enumerate}
\item If $\tau_2^\infty>\tau_1^\infty$, then $\sigma_1^\infty=\tau_1^\infty$ and the claim is proved because the corresponding stage $u=0$  does not occur. Indeed, if by contradiction $(\sigma_1^\infty,\tau_1^\infty)\ne\varnothing$ then, up to a subsequence, we would have 
    $(\sigma_1^{T_n},\tau_1^{T_n})\ne\varnothing$ for every $n$ large enough. But, looking at Theorem \ref{B0soc}, this 
    implies  $\tau_1^{T_n}=\tau_2^{T_n}$ and hence $\tau_2^\infty=\tau_1^\infty$ against our assumption. 
\item If $\tau_2^\infty=\tau_1^\infty$ hen $u^\infty=0-u_{\max}-0$ and the claim is true with $\tau_1=\tau_2$.
\end{enumerate}
Then the optimal control is always of the form \eqref{ocIHBf} and the theorem is completely proved. \QED

\section{Bocop simulations}\label{SectBocop}

To conclude, we present some numerical simulations done by using the Bocop package, \cite{Bocop, BocopExamples}.
For simplicity, they are made on a finite time interval $[0,T]$ where $T$ is taken to be large enough to ensure that in the last part of the epidemic horizon the herd immunity is reached, that is, $i$ is decreasing and the optimal control is $0$.

The Bocop package implements a local optimization method. The optimal control problem is approximated by a finite dimensional optimization problem (NLP) using a time discretization. The NLP problem is solved by the well known software Ipopt, using sparse exact derivatives computed by CppAD.
From the list of discretization formulas proposed by the package, 
we always obtained good results by using Gauss II implicit methods with 4000 time steps.

To the reader who wants to reproduce our simulations, we suggest to create a starting point by running the solver a first time without imposing the state constraint (i.e.\ the upper bound on $i$) and for $\lambda_2=0$. This starting point must then be used in solving the problem with full constraints and gradually increasing values of $\lambda_2$ by selecting the option "Re-use solution from .sol file below, including multipliers" in the Optimization menu. 
The reader is warned that neglecting this suggestion could result in a frustrating series of stops due to local infeasibililty. 
Another strategy to overcome local infeasibility errors is to enlarge the time horizon, till to several times the time expected to reach herd immunity. 

Our simulations fall in the framework of Example \ref{ex_covita} about Covid-19 epidemic in Italy and provide numerical solutions to the optimal control problem $\mathcal{P}^\infty_{{\lambda_1},{\lambda_2}}$ (see Problem \ref{IHProb})
in different times of the epidemic and different ratios $\lambda_2/\lambda_1$. In fact, we fix $\lambda_1=1$ and take several different values of $\lambda_2$.  The first two example deal with the first part of the epidemic with and without the ICU constraint. The third is about the second part of the epidemic.

\begin{figure}[H]
\begin{center}
\includegraphics[width=0.5\textwidth]{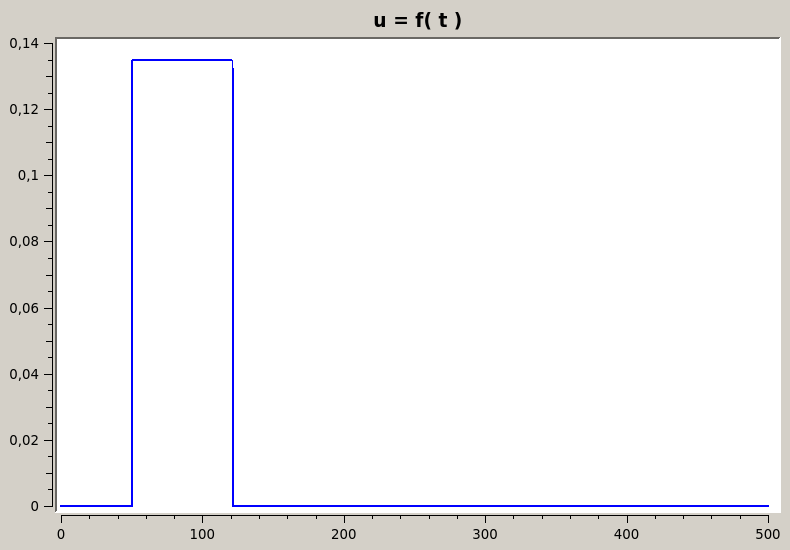}\\
\includegraphics[width=0.5\textwidth]{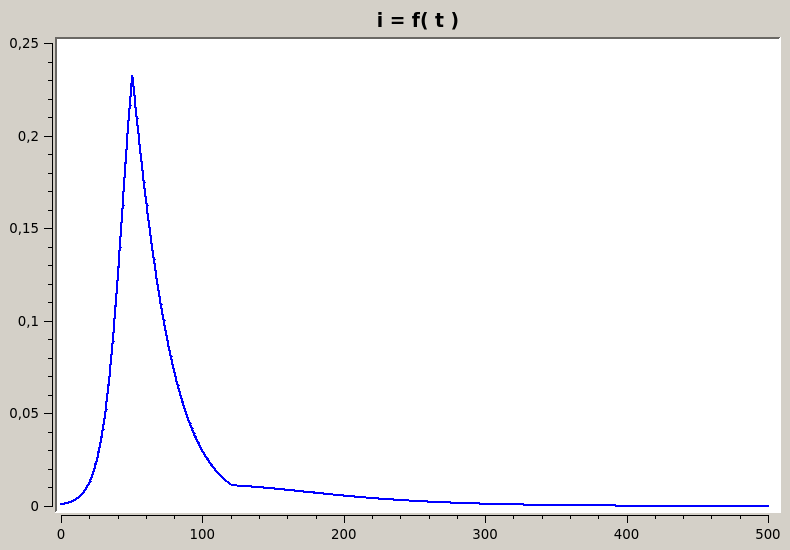}
\end{center}
\caption{Optimal control and state of infections for $\lambda_2=6$}
\label{Figure_Abang}
\end{figure}

\begin{example}{\em Our first simulation concerns the initial part of the Covid-19 epidemic in Italy 
with  epidemic parameters $\gamma=0.0714$, $\beta=0.2142$, $u_{\max}=0.135$  (see Example \ref{ex_covita}) without imposing any ICU constraint. 
We consider a time horizon of 500 days. Computations are made starting from $s_0=0.94$ and $i_0=0.001$. 

A numerical simulation with $\lambda_2=6$ 
is displayed in Figure \ref{Figure_Abang}. With the same data but $\lambda_2\in[0,2]$ the optimal control turns out to be $u\equiv0$. 
This numerical simulation illustrates the results stated in Theorem \ref{ThExtrA}, Corollary \ref{corThExtrA},  Theorem \ref{remIM=1} and Theorem \ref{IHA}.
}
\end{example}

\begin{figure}[H]
\begin{center}
\includegraphics[width=0.5\textwidth]{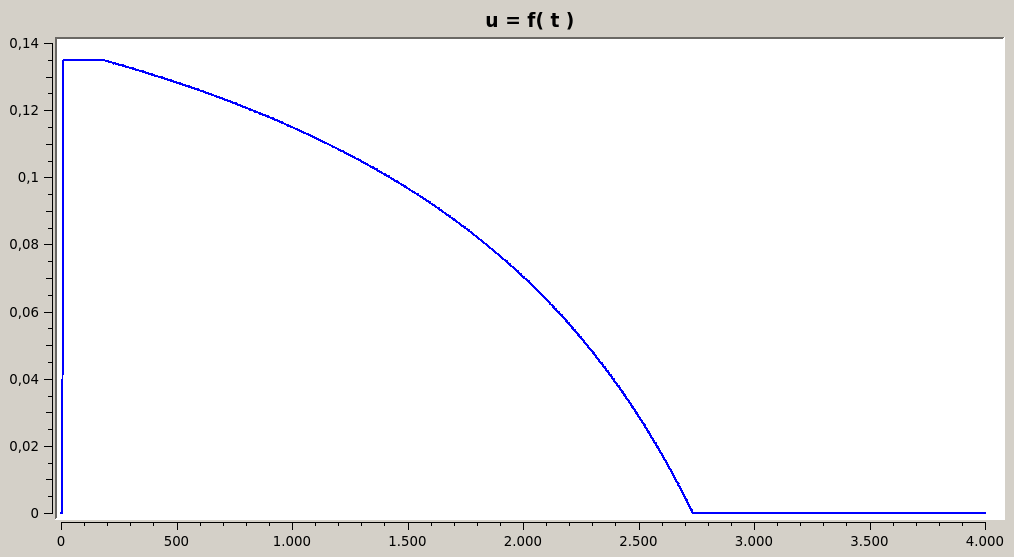}\\
\includegraphics[width=0.5\textwidth]{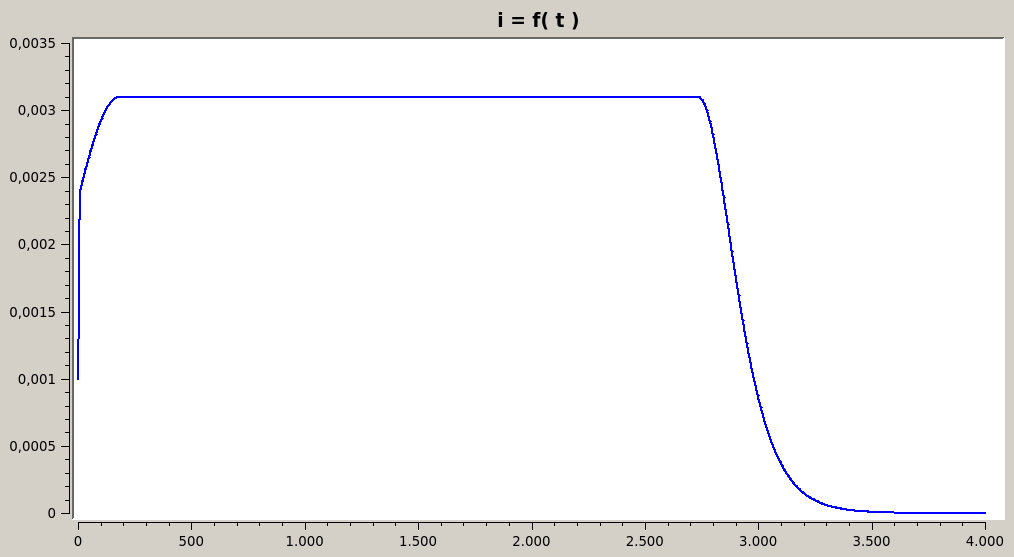}
\end{center}
\caption{Optimal control and state of infections for $\lambda_2\in[0,8]$}
\label{l1=0}
\end{figure}

\begin{example}\label{ita_cov_noICU}{\em Our second simulation differs from the first simply because we now impose the ICU constraint with  
$i_M=0.0031$. The viability regions and the corresponding  position of the initial point $(s_0,i_0)=(0.94,0.001)$  are shown in Figure \ref{ita_ini}. 
For $\lambda_2\in[0,8]$, 
the solver converged always to the greedy solution (displayed in the Figure \ref{l1=0}).
For a greater choice of $\lambda_2$ we have not been able to obtain a numerical solution.
}
\end{example}

\begin{figure}[H]
\begin{center}
\includegraphics[width=0.5\textwidth]{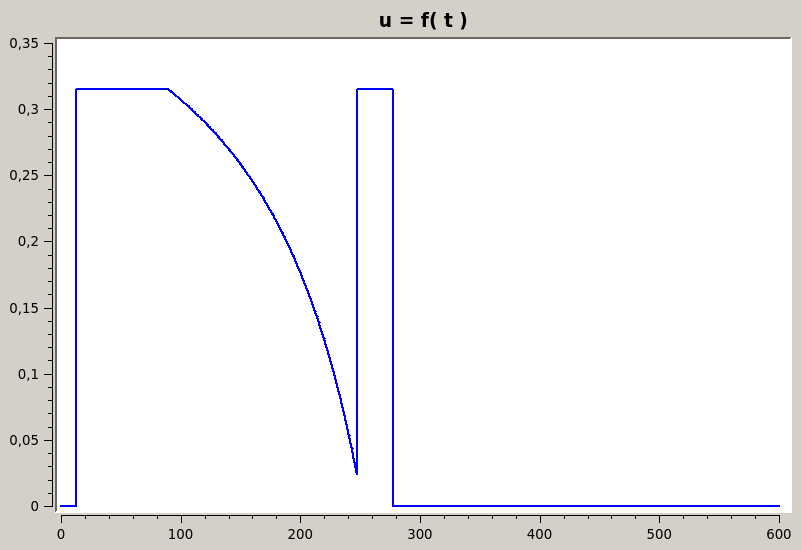}\\
\includegraphics[width=0.5\textwidth]{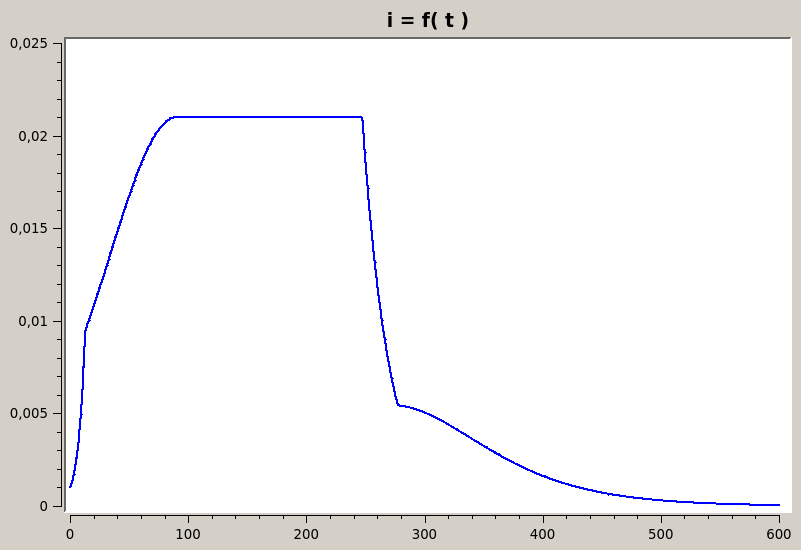}
\end{center}
\caption{Optimal control and state of infections for $\lambda_2=30$}
\label{ITA1}
\end{figure}

\begin{example}[Italy]\label{it2}{\em Our last simulation concerns the Covid-19 epidemic state of autumn 2021 in Italy considered in the second part of Example \ref{ex_covita}. 
The epidemic parameters are $\gamma=0.0714$, $\beta=0.5$, $i_M=0.021$, $u_{\max}=0.315$. 
Computations are done starting from $s_0=0.5$ and $i_0=0.001$. The viability regions and the corresponding  position of the initial point $(s_0,i_0)$ are shown in Figure \ref{ita_aut}.
We consider a time horizon of 600 days. 
For $\lambda_2\in[0,19]$ the solver converged always to the greedy solution. From $\lambda_2=20$ a second $u_{\max}$ regime starts to appear.  The graphs of the optimal control and the corresponding state of infections are displayed in Figure \ref{ITA1} for $\lambda_2=30$. 
}
\end{example}

\section{Conclusions}

We studied the infinite horizon optimal control problem for a SIR epidemic with an ICU contraint. 
We proved that, starting from an initial epidemic state $(s_0,i_0)$ in the viable set, there exists an optimal control of the form
\begin{equation*}
u(t)=\begin{cases}
0&\mbox{ if }t\in(0,\tau_0),\\
\umax&\mbox{ if }t\in(\tau_0,\tau_1),\\
\dis\beta-\frac{\gamma}{s(\tau_2)+\gamma i_M(\tau_2-t)}&\mbox{ if }t\in(\tau_1,\tau_2),\\[2ex]
\umax&\mbox{ in }(\tau_2,\tau_3),\\
0&\mbox{ in }(\tau_3,\infty), 
\end{cases}
\end{equation*}
for a suitable choice of times $0\le\tau_0\le \tau_1\le\tau_2\le\tau_3<\infty$. 
Some numerical simulations have been performed to illustrate this result. Finally, let us emphasize that the problem under consideration in which purely confinement-related costs $\lambda_1 u$ are completed with social security (convalescence) costs $\lambda_2 i$ provides a multi strict confinement scenario (illustrated by the presence of the two $u_{\max}$ branches).

\bibliographystyle{abbrv}

\bibliography{Pare37}

\end{document}